\documentclass[12pt]{article}
\usepackage{diagbox}
\usepackage{amssymb}
\usepackage{amsmath,bm}
\usepackage{multirow}
\usepackage{makecell}
\usepackage{graphicx}
\usepackage{subfigure}
\usepackage{float}
\usepackage[titletoc, title]{appendix}
\usepackage{stmaryrd}
\usepackage{lipsum}
\usepackage{amsfonts}
\usepackage{amsthm}
\usepackage{graphicx}
\usepackage{epstopdf}
\usepackage{algorithmic}
\usepackage{amsopn}
\usepackage[nameinlink]{cleveref}
\usepackage{amssymb,amsmath,amscd}
\usepackage[margin=1in]{geometry}
\usepackage{color}

\newtheorem{corollary}{Corollary}[section]
\newtheorem{remark}{Remark}[section]
\newtheorem{theorem}{Theorem}[section]
\newtheorem{lemma}{Lemma}[section]

\begin{document}

\title{Analysis of a hybridizable discontinuous Galerkin  method for the Maxwell operator}
\date{}
\author{
	Gang Chen\thanks{School of Mathematics Sciences, University of Electronic Science and Technology of China, Chengdu 611731, China. email: cglwdm@uestc.edu.cn. The first author's research is supported by
		\color{black} National Natural Science Foundation of China (NSFC) grant no. 11801063, China Postdoctoral Science Foundation grant no. 2018M633339, and Key Laboratory of Numerical Simulation of Sichuan Province (Neijiang, Sichuan Province) grant no. 2017KF003.}, 
	Jintao Cui\thanks{Department of Applied Mathematics, The Hong Kong Polytechnic University,
		Hung Hom, Hong Kong.
		email:jintao.cui@polyu.edu.hk.
		The second author's research is supported in part by \color{black} the Hong Kong RGC, General Research Fund (GRF) grant no. 15302518 and the National Natural Science Foundation of China (NSFC) grant no. 11771367.}, 
	Liwei Xu\thanks{School of Mathematics Sciences,
		University of Electronic Science and Technology of China,
		Chengdu 611731, China.  Corresponding author: email:xul@uestc.edu.cn. The third author's research is supported in part by a Key Project of the Major Research Plan of NSFC grant no. 91630205 and the NSFC grant no. 11771068.}
	}
\maketitle
% REQUIRED
\begin{abstract} 
In this paper, we study a hybridizable discontinuous Galerkin (HDG) method for the Maxwell operator. The only global unknowns are defined on the inter-element boundaries, and the numerical solutions are obtained by using discontinuous polynomial approximations. The error analysis is based on a mixed curl-curl formulation for the Maxwell equations. Theoretical results are obtained under a more general regularity requirement. In particular for the low regularity case, special treatment is applied to approximate data on the boundary. The HDG method is shown to be stable and convergence in an optimal order for both high and low regularity cases. Numerical experiments with both smooth and singular analytical solutions are performed to verify the theoretical results.

keywords: Maxwell equations, HDG method, low regularity
\end{abstract}

\section{Introduction}\label{}
Let $\Omega$ be a bounded simply-connected Lipschitz polyhedron in $\mathbb{R}^3$ with connected boundary $\color{black}\Gamma:=\partial\Omega$. We consider the following Maxwell model problem:
find the vector field $\bm{u}$ and \color{black} the Lagrange multiplier \color{black} $p$ such that
\begin{eqnarray} \label{o1}
\left\{
\begin{aligned}
\nabla\times\nabla\times \bm{u} +\nabla p&=\bm{f}\,\,\,\,\,\text{ in }\Omega, \label{source}\\
\nabla\cdot\bm{u}&=g\,\,\,\,\,\,  \text{ in }\Omega,  \\
\bm{n} \times \bm{u} &=\bm{g}_T\,\, \text{ on }\Gamma,\\
p&=0\,\, \,\,\,\,\text{ on }\Gamma.
\end{aligned}
\right.
\end{eqnarray}
Here, $\bm{n} $ is the outward normal unit vector to its boundary $\Gamma$, $\bm{f}\in [L^2(\Omega)]^3$ is an external source filed, %%and $\nabla\cdot\bm{f}=0$, $g\in L^2(\Omega)$ and
$\bm{g}_T\in \bm{H}^{-\frac{1}{2}}({\color{black}\rm div}_{\tau};\Gamma)\cap[H^{\delta}(\Gamma)]^3$ %%\bm{H}^{-\frac{1}{2}}({\color{black}\rm div}_{\tau};\Gamma)$
is a given function with $\delta\in (0,1/2]$, and $\bm{H}^{-\frac{1}{2}}({\color{black}\rm div}_{\tau};\Gamma)$ is the range space of ``tangential trace" of space $\bm{H}({\color{black}\rm curl};\Omega)$, see \cite{Buffa} for the detailed description of space $\bm{H}^{-\frac{1}{2}}({\color{black}\rm div}_{\tau};\Gamma)$.

The computational electromagnetics is of great importance in many areas of engineering and science, such as aerospace industry, telecommunication, medicine, biology, etc.
Many computational techniques  have been developed for solving Maxwell equations in both frequency and time domains. The main difficulties and challenges in designing accurate, robust and efficient numerical methods for Maxwell equations of the formulation \eqref{source} are as follows.
\begin{itemize}
	\item The irregular geometries in the real-word electromagnetic problems may lead to low regularity of true solutions, which makes the designing of stable and accurate numerical scheme and its error analysis more complicated.
	%\item When $g=0$, the incompressible condition $\nabla \cdot\bm{u}=0$ resembles an important physical property. Therefore, numerical schemes are desired to preserve this condition (at least locally).
	\item The using of Lagrange multiplier in the formulation requires an inf-sup condition in numerical analysis, in order to ensure stability and uniqueness of the approximation to $p$.
	\item The presence of nonzero $\bm{g}_T$ brings difficulties into the derivation of regularity results and error analysis. Moreover, it also affects the numerical performance of the schemes.
\end{itemize}

The finite element method is one of the most popular computational techniques for solving Maxwell equations due to its great advantages in handling complex geometries.
The methods using $\bm{H}({\color{black}\rm curl})$-conforming edge elements have been studied in vast literatures, see \cite{Ned1,Ned2,MR2009375,Monk,Jin} for examples.

Since the late 1970s, discontinuous Galerkin (DG) methods have become increasingly popular in scientific computation due to their attractive features. The first DG method for solving the time-harmonic Maxwell equations in electric field was studied in \cite{LDG-Maxwell}, where a local discontinuous Galerkin (LDG) scheme was proposed for the low frequency problem. An interior penalty discontinuous Galerkin (IPDG) method for Maxwell equations with high frequency was studied in \cite{IPDG-Maxwell}, where the constants in the stability estimates and the error estimates were highly dependent on the frequency. Mixed DG approach for time-harmonic Maxwell equations with low frequency was studied in  \cite{MDG-Maxwell}. Numerical analysis has been carried out for both cases when the analytical solution is smooth and holds minimal regularity, such that
\begin{align}
\bm{u}\in [H^{s}(\Omega)]^3,\quad \nabla\times\bm{u}\in [H^{s}(\Omega)]^3,\quad p\in H^{1+s}(\Omega),\quad\bm{g}_T\in [H^{s+\frac{1}{2}}(\Gamma)]^3, \label{minimal}
\end{align}
with $s> 1/2$. Note that in model problem \eqref{source}, the domain is assumed to be a bounded simply-connected Lipschitz polyhedron, and the relative magnetic permeability and the relative electric permittivity are set to be one; hence the above regularity results can be achieved by setting $\bm{f}\in [L^2(\Omega)]^3$, $g\in L^2(\Omega)$ and $\bm{g}_T=\bm{0}$.

In recent years, a ``new" type of DG method, the hybridizable discontinuous Galerkin (HDG) method has been successfully applied to solve various types of differential equations. Compared with conventional DG methods, the HDG method has several distinct advantages: static condensation, less degree of freedom (reduced discrete system size), flexible in meshing (inherent from DG methods), easy to design and implement, local conservation of physical quantities, and so on. The first work on HDG method for time-harmonic Maxwell equations was founded in \cite{Cockburn-Maxwell}. Two types of HDG schemes were introduced therein and numerical experiments were presented to show the performance of proposed schemes. At a later time, HDG methods for the time-harmonic Maxwell equations with zero frequency were introduced and analyzed in \cite{Shike-Maxwell}. The regularity assumptions for the convergence analysis are much more restrictive than \eqref{minimal}, i.e.,
\begin{align*}
\bm{u}\in [H^{t}(\Omega)]^3, \quad\nabla\times\bm{u}\in [H^{s}(\Omega)]^3, \quad p\in H^r(\Omega),
\end{align*}
with $t,s,r\ge 1$. This assumption is too strong to be satisfied (hence not realistic in practice). Actually, it can not even be met for problems defined on bounded simply-connected Lipschitz polyhedral domains. \color{black}  Later in \cite{MR3771897}, an HDG method was proved to converge optimally in both energy norm and $L^2$-norm for $\bm u$ under the regularity \eqref{minimal}, and {\em a posteriori} error estimates were derived.
\color{black}

Recently in \cite{IBC1,IBC2}, the authors studied the HDG methods for time-harmonic Maxwell equations with {\em impedance boundary conditions} and {\em large wavenumber}. The error analysis in both work was derived by using the regularity results given in \cite{Maxwell-IBC,MR3265182}, where the constants in error estimates were shown to be explicitly dependent on the wave number.  We also refer to \cite{MR2407321} for a nonconforming finite element method solving the two-dimensional curl-curl (Maxwell operator)  problem,
refer to \cite{MR3614286} for a hybrid-mesh HDG method to the time-harmonic Maxwell’s equations, 
refer to \cite{MR2579978,MR2966917} for preconditioners for the discretized
time-harmonic Maxwell's equations in mixed form,
and refer to \cite{MR2737418} for an initial-boundary value problem for the Maxwell's equations.

\color{black}
The aim of this paper is to propose and analyze a new HDG method for solving the Maxwell operator. 
We first derive the regularity results of Maxwell's equations with general given data $\bm{f}$, $g$ and $\bm{g}_T$, which are new in the literature: there exists a constant $s\in (1/2,1/2+\delta]$ $(\delta>0)$ such that
\begin{align*}
\|\bm{u}\|_{s}
\le C
\left( \|\bm{f}\|_0+\|g\|_0+\|\bm{g}_T\|_{\delta,\Gamma} +\|\bm{g}_T\|_{-\frac{1}{2},{\rm div}_{\tau}}
\right),
\end{align*}
and
\begin{align*}
\|\nabla\times\bm{u}\|_{s}\le C 
\left( \|\bm{f}\|_0+\|g\|_0+\|\bm{g}_T\|_{1+\delta,\Gamma} +\|\bm{g}_T\|_{-\frac{1}{2},{\rm div}_{\tau}}
\right).
\end{align*}
Then we introduce a new HDG method based on particularly designed projections on the boundary data, with which the oscillation terms corresponding to the boundary data can be eliminated.
We show that the HDG method is stable and has optimal error convergence rates in mesh-dependent energy norm, under a general (more realistic) regularity requirements. The general regularity is the lowest in the existing literature. Move precisely, we only require that
\begin{align*}
\bm{u}\in [H^{s}(\Omega)]^3,\quad
\nabla\times\bm{u}\in [H^{s}(\Omega)]^3,\quad
p\in H^{1+s}(\Omega),\\
\bm{g}_T\in \bm{H}^{-\frac{1}{2}}({\rm div}_{\tau};\Gamma)\cap[H^{\delta}(\Gamma)]^3,\quad s>1/2,
\quad
\delta >0.
\end{align*}

\color{black}
The rest of this paper is organized as follows. In section 1, we introduce the notation and give the regularity results. In section 2, we propose an HDG method for the Maxwell model problem \eqref{source}. The stability analysis and convergence analysis for the proposed HDG scheme are given in sections 3 and 4, respectively. In section 5, we perform some numerical experiments to verify the theoretical results. Finally, we conclude this paper in section 6.

Throughout this paper, we use $C$ to denote a positive constant independent of mesh size, which may take on different values at each occurrence. We also use $a\sim b$ to stand for $Cb\le a\le Cb$.

\section{Regularity analysis}

For any bounded domain $\Lambda\color{black}\subset\color{black} \mathbb{R}^s$ $(s=1,2,3)$, let $H^{m}(\Lambda)$ denote the usual  $m^{th}$-order Sobolev space on $\Lambda$, and $\|\cdot\|_{m, \Lambda}$, $|\cdot|_{m,\Lambda}$  denote the norm and semi-norm on $H^{m}(\Lambda)$.
We use $(\cdot,\cdot)_{m,\Lambda}$ to denote the inner product of $H^m(\Lambda)$, with $(\cdot,\cdot)_{\Lambda}:=(\cdot,\cdot)_{0,\Lambda}$.
When $\Lambda=\Omega$, we denote $\|\cdot\|_{m }:=\|\cdot\|_{m, \Omega}, |\cdot|_{m}:=|\cdot|_{m,\Omega}$, $(\cdot,\cdot):=(\cdot,\cdot)_{\Omega}$.
In particular, when $\Lambda\in \mathbb{R}^{2}$, we use $\langle\cdot,\cdot\rangle_{\Lambda}$ to replace $(\cdot,\cdot)_{\Lambda}$;
when $\Lambda\in \mathbb{R}^{1}$, we use $\langle\!\langle\cdot,\cdot\rangle\!\rangle_{\Lambda}$ to replace $(\cdot,\cdot)_{\Lambda}$.
Note that bold face fonts will be used for vector (or tensor) analogues of the Sobolev spaces along with vector-valued (or tensor-valued) functions.
We denote the following Sobolev spaces
\begin{align*}
\bm{H}({\color{black}\rm curl};\Omega)&:=\{\bm{v}\in [L^2(\Omega)]^3: \nabla\times\bm{v}\in [L^2(\Omega)]^3 \},\\
\bm{H}^s({\color{black}\rm curl};\Omega)&:=\{\bm{v}\in [H^s(\Omega)]^3: \nabla\times\bm{v}\in [H^s(\Omega)]^3 \}\,\text{ with }s\ge 0,\\
\bm{H}_0({\color{black}\rm curl};\Omega)&:=\{\bm{v}\in \bm{H}({\color{black}\rm curl};\Omega): \bm{n}\times\bm{v}=\bm{0} \,\,\text{on}\,\,\Gamma\},\\
\bm{H}({\color{black}\rm div} ;\Omega)&:=\{\bm{v}\in [L^2(\Omega)]^3: \nabla\cdot\bm{v}\in L^2(\Omega) \},\\
\bm{H}_0({\color{black}\rm div} ;\Omega)&:=\{\bm{v}\in \bm{H}({\color{black}\rm div} ;\Omega): \bm{n}\cdot\bm{v}=0 \,\,\text{on}\,\,\Gamma \},\\
\bm{H}({\color{black}\rm div} ^0;\Omega)&:=\{\bm{v}\in \bm{H}({\color{black}\rm div} ;\Omega): \nabla\cdot\bm{v}=0\},
\end{align*}
and
\begin{align*}
\bm{X} :=\bm{H}({\color{black}\rm curl};\Omega)\cap \bm{H}({\color{black}\rm div} ;\Omega),\qquad
\bm{X} _N:=\bm{H}_0({\color{black}\rm curl};\Omega)\cap \bm{H}({\color{black}\rm div} ;\Omega),\\
%%\bm{X} _{N,0}&:=&\bm{H}_0({\color{black}\rm curl};\Omega)\bigcap \bm{H}({\color{black}\rm div}^0 ;\Omega),\\
\bm{X} _T:=\bm{H}({\color{black}\rm curl};\Omega)\cap \bm{H}_0({\color{black}\rm div} ;\Omega).
\end{align*}

We define the norms on $\bm{H}^s({\color{black}\rm curl};\Omega)$ and $\bm{H}({\color{black}\rm div}_{\tau};\Gamma)$, respectively as follows:
\begin{align*}
\|\bm{v}\|_{s,{\color{black}\rm curl}}=\left(  \|\bm{v}\|_s^2+\|\nabla\times\bm{v}\|_s^2    \right)^{\frac{1}{2}},\qquad
\|\bm{w}\|_{-\frac{1}{2},{\color{black}\rm div}_{\tau}}=\inf_{^{\bm{v}\in \bm{H}({\color{black}\rm curl};\Omega):}_{\bm{n} \times\bm{v}|_{\Gamma}=\bm{w}}} \|\bm{v}\|_{0,{\color{black}\rm curl}}.
\end{align*}

\subsection{ Regularity}

\begin{lemma} [Regularity for $p$] 
	Let $\bm f\in \bm H({\rm div};\Omega)$, then
	the model problem \eqref{source} has a unique solution $p$, and there exists a constant $s_2\in (1/2,1]$ such that
	\begin{align}\label{219}
	\|\nabla p\|_0\le \|\bm f\|_0,\qquad
	\|p\|_{1+s_2}\color{black}\le C \color{black} \|\color{black}\nabla\cdot\bm{f}\|_0.
	\end{align}
\end{lemma}
\begin{proof}
	\color{black}
	We apply $\nabla\cdot$ on the first equation of \eqref{source}, and combine the last equation of \eqref{source} to get
	\begin{subequations}\label{elliptic}
		\begin{align}
		\Delta p&=\nabla\cdot\bm f&\text{in }\Omega,\\
		p&=0&\text{on }\Gamma.
		\end{align}
	\end{subequations}
	By using integration by parts and \eqref{elliptic}, we get
	\begin{align*}
	(\nabla p,\nabla p)=-(\Delta p, p)=-(\nabla\cdot\bm f, p)=(\bm f,\nabla p),
	\end{align*}
	which implies $\|\nabla p\|_0\le \|\bm f\|_0$.
	Utilizing \cite[Corollary 2.6.7]{Grisvard} for the elliptic problem \eqref{elliptic},
	\color{black}
	there exists a unique $p$ and a regularity index $s_2\in (1/2,1]$ such that
	\begin{align*}
	\|p\|_{1+s_2}\color{black}\le C \color{black} \ \color{black} \|\color{black}\nabla\cdot\bm{f}\|_0,
	\end{align*}
	which proves the lemma.
\end{proof}

\begin{lemma}[cf. \cite{Maxwell-1997}]\label{con}
	If $\Omega$ is a Lipschitz domain, let $\bm{v}\in \bm{X} $, $\bm{n}\cdot\bm{v}\in L^2(\Gamma_{\tau})$ and $\bm{n}\times\bm{v}\in L^2(\Gamma_{\nu})$, where  $\Gamma_{\tau}\bigcap\Gamma_{\nu}=\emptyset$, $\Gamma_{\tau}\bigcup\Gamma_{\nu}=\Gamma$, then we have
	\begin{align*}
	\|\bm{v}\|_0\color{black}\le C \color{black}
	\left(
	\|\nabla\times\bm{v}\|_0+\|\nabla\cdot\bm{v}\|_0
	+\|\bm{n}\times\bm{v}\|_{0,\Gamma_{\tau}}+\|\bm{n}\cdot\bm{v}\|_{0,\Gamma_{\nu}}
	\right).
	\end{align*}
\end{lemma}

\begin{lemma} The Maxwell model problem \eqref{source} has a unique solution $\bm{u}$, and the following stability holds:
	\begin{align*}
	\|\bm{u}\|_{0,{\color{black}\rm curl}}+
	\|\nabla\times\bm{u}\|_{0,{\color{black}\rm curl}}
	\color{black}\le C \color{black}
	\left(
	\|\bm{f}\|_0+\|g\|_0
	+\|\bm{g}_T\|_{-\frac{1}{2},{\color{black}\rm div}_{\tau}}+\|\bm{g}_T\|_{0,\Gamma}
	\right).
	\end{align*}
\end{lemma}
\begin{proof}
	A natural variational problem of \eqref{source} reads: for all $(\bm{v},q)\in \bm{H}_0({\color{black}\rm curl};\Omega)\times H^1_0(\Omega)$, find $(\bm{u},q)\in \bm{H}({\color{black}\rm curl};\Omega)\times H^1_0(\Omega)$ such that
	\begin{subequations}
		\begin{align}
		(\nabla\times\bm{u},\nabla\times\bm{v})+(\nabla p,\bm{v})&=(\bm{f},\bm{v}),\label{var1}\\
		-(\bm{u},\nabla q)&=(g,q),\\
		\bm{n} \times \bm{u}|_{\Gamma}&=\bm{g}_T.\label{var3}
		\end{align}
	\end{subequations}
	Since $\bm{g}_T\in\bm{H}^{-\frac{1}{2}}({\color{black}\rm div}_{\tau};\Gamma)$, and $\bm{H}^{-\frac{1}{2}}({\color{black}\rm div}_{\tau};\Gamma)$ is the range space of ``tangential trace" of space $\bm{H}({\color{black}\rm curl};\Omega)$, there exists a $\bm{u}_0\in \bm{H}({\color{black}\rm curl};\Omega)$ such that
	\begin{eqnarray}\label{uoest}
	\bm{n} \times\bm{u}_0|_{\Gamma}=\bm{g}_T.
	\end{eqnarray}
	Setting $\bm{w}=\bm{u}-\bm{u}_0$, the variational problem \eqref{var1}--\eqref{var3} becomes:
	for all $(\bm{v},q)\in \bm{H}_0({\color{black}\rm curl};\Omega)\times H^1_0(\Omega)$, find $(\bm{w},q)\in \bm{H}_0({\color{black}\rm curl};\Omega)\times H^1_0(\Omega)$ such that
	\begin{subequations}
		\begin{align}
		(\nabla\times\bm{w},\nabla\times\bm{v})+(\nabla p,\bm{v})&=(\bm{f},\bm{v})-(\nabla\times\bm{u}_0,\nabla\times\bm{v}),\label{varv1}\\
		-(\bm{w},\nabla q)&=(g,q)+(\bm{u}_0,\nabla q)\label{varv2}.
		\end{align}
	\end{subequations}
	Taking $(\bm{v},q)=(\bm{w},p)\in \bm{H}_0({\color{black}\rm curl};\Omega)\times H^1_0(\Omega)$ in \eqref{varv1}--\eqref{varv2} and adding them together yield
	\begin{align}\label{uopest1}
	\|\nabla\times\bm{w}\|^2_0&=(\bm{f},\bm{w})-(\nabla\times\bm{u}_0,\nabla\times\bm{w})+(g,p)+(\bm{u}_0,\nabla p)\nonumber\\
	&\le C_{\beta}(\|\bm{f}\|_0+\|g\|_0+\|p\|_1+\|\bm{u}_0\|_{0,{\color{black}\rm curl}})^2+\frac{\beta^2}{4}\|\bm{w}\|^2_0+\frac{3}{4}\|\nabla\times\bm{w}\|^2_0,
	\end{align}
	where $\beta>0$ is a constant  to be specified later, and $C_{\beta}>0$ is a constant dependent on $\beta$.
	It then follows from \eqref{219}, \eqref{uoest}, \eqref{uopest1} and the definition of $\|\cdot\|_{-\frac{1}{2},{\color{black}\rm div}_{\tau}}$ that
	\begin{eqnarray} \label{232}
	\|\nabla\times\bm{w}\|_0\le C_{\beta}(\|\bm{f}\|_0+\|g\|_0
	+\|\bm{g}_T\|_{-\frac{1}{2},{\color{black}\rm div}_{\tau}})+\beta\|\bm{w}\|_0.
	\end{eqnarray}
	By using a triangle inequality and Lemma \ref{con}, we can get
	\begin{align}\label{west1}
	\|\bm{w}\|_0&\le \|\bm{u}\|_0+ \|\bm{u}_0\|_0 \nonumber\\
	&\le C( \|\nabla\times\bm{u}\|_0 +\|g\|_0+\|\bm{g}_T\|_{0,\Gamma}  )  +\|\bm{u}_0\|_0 \nonumber\\
	&\le C_0( \|\nabla\times\bm{w}\|_0+\|\nabla\times\bm{u}_0\|_0 +\|g\|_0+\|\bm{g}_T\|_{0,\Gamma}  )  +\|\bm{u}_0\|_0.
	\end{align}
	Combining \eqref{uoest}, \eqref{232}, \eqref{west1} and the definition of $\|\cdot\|_{-\frac{1}{2},{\color{black}\rm div}_{\tau}}$, we get
	\begin{align*}
	\|\bm{w}\|_0 \le C_{\beta}(\|\bm{f}\|_0+\|g\|_0+\|\bm{g}_T\|_{-\frac{1}{2},{\color{black}\rm div}_{\tau}}+\|\bm{g}_T\|_{0,\Gamma})+C_0\beta\|\bm{w}\|_0.
	\end{align*}
	Hence by taking $\beta=\frac{1}{2C_0}$, we arrive at
	\begin{align*}
	\|\bm{w}\|_0\color{black}\le C
	\left( \color{black}\|\bm{f}\|_0+\|g\|_0
	+\|\bm{g}_T\|_{-\frac{1}{2},{\color{black}\rm div}_{\tau}}+\|\bm{g}_T\|_{0,\Gamma}\right).
	\end{align*}
	Then it directly follows from a triangle inequality that
	\begin{align*}
	\|\bm{u}\|_{0,{\color{black}\rm curl}}&\le \|\bm{w}\|_{0,{\color{black}\rm curl}}+\|\bm{u}_0\|_{0,{\color{black}\rm curl}}\\
	&\le C \color{black}
	\left( \|\bm{f}\|_0+\|g\|_0+\|\bm{g}_T\|_{-\frac{1}{2},{\color{black}\rm div}_{\tau}}+\|\bm{g}_T\|_{0,\Gamma}\right)+\|\bm u_0\|_{0,\rm curl}.
	\end{align*}
	Using the definition of $\|\cdot\|_{-\frac{1}{2},{\color{black}\rm div}_{\tau}}$ and the arbitrary feature of $\bm u_0$, we get
	\begin{align*}
	\|\bm{u}\|_{0,{\color{black}\rm curl}}\le C \color{black}
	\left( \|\bm{f}\|_0+\|g\|_0+\|\bm{g}_T\|_{-\frac{1}{2},{\color{black}\rm div}_{\tau}}+\|\bm{g}_T\|_{0,\Gamma}\right).
	\end{align*}
	Since $\|\nabla\times\nabla\times\bm{u}\|_0=\|\bm{f}-\nabla p\|_0\color{black}\le 2 \color{black}\|\bm{f}\|_0$, we have
	\begin{align*}
	\|\nabla\times\bm{u}\|_{0,{\color{black}\rm curl}}\color{black}\le C \color{black}
	\left(\|\bm{f}\|_0+\|g\|_0 +\|\bm{g}_T\|_{-\frac{1}{2},{\color{black}\rm div}_{\tau}}+\|\bm{g}_T\|_{0,\Gamma}\right),
	\end{align*}
	which completes the proof.
\end{proof}

The following Sobolev embedding results holds true.
\begin{lemma}[cf. \cite{Vector}] \label{embed}
	If $\Omega$ is a Lipschitz polyhedral domain, there exists a real number ${\color{black} \sigma}\in (1/2,1]$ such that $\bm{X}_T$ and $\bm{X}_N$ are continuously imbedded in $[H^{{\color{black} \sigma}}(\Omega)]^3$.
\end{lemma}

\begin{lemma}[cf. \cite{Regularity}] \label{embed2}
	If $\Omega$ is a Lipschitz polyhedral domain, then for each $\delta\in (0,1/2)$, the spaces
	\begin{align*}
	\bm{W}:=\{\bm{w}\in \bm{X} :  (\bm{n} \times\bm{w})|_{\Gamma}\in [H^{\delta}(\Gamma)]^3\}\,
	\text{ and }\,\,
	\bm{V}:=\{\bm{w}\in \bm{X} :  (\bm{n} \cdot\bm{w})|_{\Gamma}\in [H^{\delta}(\Gamma)]^3\}
	\end{align*}
	are continuously imbedded in $[H^{s_1}(\Omega)]^3$, where $s_1=\min(\delta+1/2,{\color{black} \sigma})$ and ${\color{black} \sigma}\in (1/2,1]$ are defined as in Lemma \ref{embed}.
\end{lemma}

Now we are well prepared to derive the following regularity result for problem \eqref{source}.
\begin{theorem} \label{reg} If $\Omega$ is a Lipschitz polyhedral domain and $\bm f\in \bm H({\rm div};\Omega)$, the model problem \eqref{source} has a unique solution $(\bm{u},p)$. Moreover, there exists a regularity index $s\in(1/2,1/2+\delta]$, such that
	\begin{align*}
	\|\bm{u}\|_{s}
	\color{black}\le C \color{black}
	\left( \|\bm{f}\|_0+\|g\|_0+\|\bm{g}_T\|_{\delta,\Gamma} +\|\bm{g}_T\|_{-\frac{1}{2},{\color{black}\rm div}_{\tau}}
	\right)
	,\qquad
	\color{black}\|p\|_{1+s}\le C\|\nabla\cdot\bm f\|_0.
	\end{align*}
	If in addition, $\bm g_T\in [H^{1+\delta}(\Omega)]^3$, there holds
	\begin{align*}
	\|\nabla\times\bm{u}\|_{s}\color{black}\le C \color{black}
	\left( \|\bm{f}\|_0+\|g\|_0+\|\bm{g}_T\|_{1+\delta,\Gamma} +\|\bm{g}_T\|_{-\frac{1}{2},{\color{black}\rm div}_{\tau}}
	\right).
	\end{align*}
\end{theorem} \color{black}
%
%%%%%%%%%%%%%%%%%%%%%%%%%%%%%%%%%%%%%%%%%%%%%%%%%%%%%%%

\section{HDG finite element method}

Let $\mathcal{T}_h=\bigcup\{T\}$ be the conforming triangulation of $\Omega$ made of shape-regular simplicial elements. For each $T\in\mathcal{T}_h$, we let $h_T$ be the infimum of the diameters of spheres containing $T$ and denote the mesh size $h:=\max_{T\in\mathcal{T}_h}h_T$. Let $\mathcal{F}_h=\bigcup\{F\}$ be the union of all  faces of $T\in\mathcal{T}_h$, and let $\mathcal{F}_h^I$ and $\mathcal{F}_h^B$ be all the interior faces and boundary faces, respectively. We denote by $h_F$ the length of the diameter of circle containing face $F$.
For all $T\in\mathcal{T}_h$ and $F\in\mathcal{F}_h$, we denote by $ \bm{n}_T $ and $\bm{n}_F$   the  unit outward  normal vectors along $\partial T$ and  face $F$, respectively. Broken curl, div and gradient operators with respect to mesh partition $\mathcal{T}_h$ are donated by $\nabla_h\times$, $\nabla_h\cdot$ and $\nabla_h$, respectively.
For $u,v\in L^2(\partial\mathcal{T}_h)$, we define the following  inner product and norm
\begin{align*}
\langle u,v \rangle_{\partial\mathcal{T}_h}=\sum_{T\in\mathcal{T}_h}\langle u,v\rangle_{\partial T},\qquad
\|v\|^2_{0,\partial\mathcal{T}_h}=\sum_{T\in\mathcal{T}_h}\|v\|^2_{0,\partial T}.
\end{align*}
\color{black}
Let $F=\partial T \cap \partial T^{'}$ be an interior face shared by element $T$ and element $T^{'}$, and $\bm n_F$ be the unit normal pointing from $T$ to $T^{'}$. For any piecewise smooth function $\bm{\phi}$, we
define the jump of $\bm{\phi}$ on $F$ as
\begin{align*}
[\![ \bm{\phi}]\!]:=\bm{\phi}|_T-\bm{\phi}^{'}|_{T^{'}}.
\end{align*}
On the boundary face $F=\partial T\cap \Gamma$, we set $[\![ \bm{\phi}]\!]:=\bm{\phi}$.

For a bounded domain $\Lambda\color{black}\subset\color{black} \mathbb{R}^s$ $(s=1,2,3)$ and for an integer $k\ge 0$, $\mathbb{P}_{k}(\Lambda)$  denotes the set of all \color{black} polynomials \color{black} defined on $\Lambda$ with the degree at most $k$. The finite element spaces $\mathbb{P}_j(\mathcal{T}_h)$ and $\mathbb{P}_j(\mathcal{F}_h)$ are defined by
\begin{align*}
\mathbb{P}_j(\mathcal{T}_h):=\{q_h\in L^2(\Omega): q_h|_T\in \mathbb{P}_j(T),\ \forall \,T\in\mathcal{T}_h\},\\
\mathbb{P}_j(\mathcal{F}_h):=\{q_h\in L^2(\mathcal{F}_h): q_h|_E\in \mathbb{P}_j(E),\ \forall\, F\in\mathcal{F}_h\}.
\end{align*}\color{black}

\subsection{Interpolations}
\subsubsection{$L^2$-projection}
For any $T\in\mathcal{T}_h$, $ F\in\mathcal{F}_h$ and any   integer $j\ge 0$, let $\Pi^{o}_j: L^2(T)\rightarrow \mathbb{P}_{j}(T)$ and $\Pi^{\partial}_j: L^2(F)\rightarrow \mathbb{P}_{j}(F)$ be the usual $L^2$ projection operators.
% from $ $ onto $\mathbb{P}_{k}(T)$, denote by $Q_{k}^b$ the $L^2$ projection from $L^2(E)$ onto $\mathbb{P}_{k}(E)$.
The following stability and approximation results are standard.
\color{black}
\begin{lemma}\label{lemma2.5} For any $T\in\mathcal{T}_h$ and $F\in\mathcal{F}_h$ \color{black} and nonnegative integer $j$\color{black}, it holds
	\begin{align*}
	\|v-\Pi^{o}_jv\|_{0,T}&\color{black}\le C \color{black} h_T^{s}|v|_{s,T}&\forall\, v\in H^{s}(T),\\
	\|v-\Pi^{o}_jv\|_{0,\partial T}&\color{black}\le C \color{black} h_T^{s-1/2}|v|_{s,T}&\forall\, v\in H^{s}(T),\\
	\|v-\Pi^{\partial}_jv\|_{0,\partial T}&\color{black}\le C \color{black} h_T^{s-1/2}|v|_{s,T}&\forall\, v\in H^{s}(T),\\
	\|\Pi^{o}_jv\|_{0,T}&\le\|v\|_{0,T}&\forall\, v\in L^{2}(T),\\
	\|\Pi^{\partial}_jv\|_{0,F}&\le\|v\|_{0,F}&\forall\, v\in L^{2}(F),
	\end{align*}
	where $s\in (1/2,j+1]$.
\end{lemma}

\subsubsection{Second type H({\color{black}\rm curl})-projection\label{curl}}
For any integer $k\ge 1$, define
\begin{align*}
\bm{\mathcal{D}}_k(T)&=[\mathbb{P}_{k-1}(T)]^3\oplus \bm{x}\mathbb{P}^*_{k-1}(T),\qquad
\bm{\mathcal{D}}_k(F)=[\mathbb{P}_{k-1}(F)]^3\oplus \bm{x}\mathbb{P}^*_{k-1}(F),
\end{align*}
where $\bm{x}=(x_1,x_2,x_3)$ and $\mathbb{P}^*_{k-1}(T)$ is the space of homogeneous polynomial of degree $k-1$. For any $T\in\mathcal{T}_h$ and face $F\subset\partial T$, we let $E\subset \partial F$ denote the edge of $T$ on face $F$. Let $\bm{t}_{E}$ be the unit tangential vector along edge $E$, and $\bm{n}_{FE}$ be the unit outer normal vector to edge $E$ on face $F$. For any $\bm{v}\in \bm{H}^s({\color{black}\rm curl};T)$ with $s>1/2$, we define the second type of $\bm H({\color{black}\rm curl})$-projection (cf. \cite{Ned2}) $\bm{\mathcal{P}}^{{\color{black}\rm curl}}_{k}: \bm{H}^s({\color{black}\rm curl};T)\to [\mathbb{P}_{k}(T)]^3$ such that
\begin{subequations}
	\begin{align*}
	\langle\!\langle\bm{\mathcal{P}}^{{\color{black}\rm curl}}_{k}\bm{v}\cdot\bm{t}_{E},w_k\rangle\!\rangle_{E}
	=\langle\!\langle\bm{v}\cdot\bm{t}_{E},w_k\rangle\!\rangle_{E}\,\,\,\,\,\, \forall\, w_k\in \mathbb{P}_k(E),\,\, \color{black}E\subset\partial F.
	\end{align*}
	For any face $F$ of $T$, there holds when $k\ge 2$,
	\begin{align*}
	\langle \bm{\mathcal{P}}^{{\color{black}\rm curl}}_{k}\bm{v},\bm{w}_{k-1}  \rangle_{F}=\langle \bm{v},\bm{w}_{k-1}  \rangle_{F}\quad\forall\, \bm{w}_{k-1}\in \bm{\mathcal{D}}_{k-1}(F),
	\end{align*}
	and when $k\ge 3$,
	\begin{align*}
	(\bm{\mathcal{P}}^{{\color{black}\rm curl}}_{k}\bm{v},\bm{w}_{k-2})=(\bm{v},\bm{w}_{k-2})\quad\forall\, \bm{w}_{k-1}\in\bm{\mathcal{D}}_{k-2}(T).
	\end{align*}
\end{subequations}
Note that the above definitions make sense when $\bm{v}\in \bm{H}^s({\color{black}\rm curl};T)$ with $s>1/2$ (see \cite[Lemma 5.1]{Regularity} for details).
\vspace{0.02in}

The following approximation properties holds true:
\begin{lemma} [cf. \cite{Ned2,Regularity,Monk}]
	For any  $\ell\in (1/2,k+1]$ and $m\in (1/2,k]$, one has
	\color{black}
	\begin{align*}
	\|\bm{v}-\bm{\mathcal{P}}^{{\color{black}\rm curl}}_{k}\bm{v}\|_{0}\color{black}&\le C  h^{\ell}\|\bm{v}\|_{\ell}&\forall \bm v\in [H^{\ell}(\Omega)]^d,\\
	\|\nabla\times(\bm{v}-\bm{\mathcal{P}}^{{\color{black}\rm curl}}_{k}\bm{v})\|_{0}\color{black}&\le C \color{black} h^{m}\|\nabla\times\bm{v}\|_{m}&\forall \bm v\in \bm H^m({\rm curl};\Omega).
	\end{align*}
	%	holds for $\bm v\in \bm H^t({\rm curl};T)$ with $t\in (1/2,k]$,
\end{lemma}
%%%%%%%
\vspace{0.02in}

\color{black}

\subsubsection{$H({\color{black}\rm curl})$-projection and $H^1$-projection on finite element spaces}

In the error analysis, we need the following $\bm H_0({\rm curl})$-conforming and $H^1$-conforming
interpolations.

\begin{lemma}[cf. {\cite[Proposition 4.5]{MR2194528}}]\label{pic}
	For any integer $k\ge 1$, let $\bm v_h\in [\mathbb P_{k}(\mathcal{T}_h)]^3$, there exists a function $\bm{\Pi}_{h,k}^{\rm curl,c}\bm v_h\in [\mathbb P_{k}(\mathcal{T}_h)]^3\cap \bm H_0({\rm curl};\Omega)
	$ such that
	\begin{align}
	\|\bm{\Pi}_{h,k}^{\rm curl,c}\bm v_h-\bm v_h\|_0&\le C \|h_F^{1/2}\bm n\times[\![\bm v_h]\!]\|_{0,\mathcal{F}_h},\label{ic1}\\
	\|\nabla_h\times(\bm{\Pi}_{h,k}^{\rm curl,c}\bm v_h-\bm v_h)\|_{0}&\le C\|h_F^{-1/2}\bm n\times[\![\bm v_h]\!]\|_{0,\mathcal{F}_h},\label{ic2}
	\end{align}
	with a constant $C>0$ independent of the mesh size.
\end{lemma}

\begin{lemma} [cf. {\cite[Theorem 2.2]{ick}}] 
	For any integer $k\ge 1$ and
	let $q_h\in \mathbb{P}_{k}(\mathcal{T}_h)$, there exists an interpolation operator \color{black} { $\mathcal{I}_{k}^c: \mathbb{P}_{k}(\mathcal{T}_h)\to \mathbb{P}_{k}(\mathcal{T}_h)\cap H^1_0(\Omega)$} \color{black}such that
	\begin{eqnarray}\label{eq:lem41}
	\|h_T^{\frac{\alpha+1}{2}}(\nabla_h q_h-\nabla \mathcal{I}_{k}^c q_h)\|_0\color{black}\le C  \|h_F^{\frac{\alpha}{2}}[\![q_h]\!]\|_{\mathcal{F}_h}.\label{ick}
	\end{eqnarray}
\end{lemma}
\begin{remark}
	The \cite[Theorem 2.2]{ick} states that \eqref{eq:lem41} holds for $\alpha=-1$. Actually its proof can be extended to any fixed
	real number $\alpha$ in a straightforward manner.
\end{remark}

\color{black}
\subsubsection{A modified projection}

Let $k\ge 1$ be an integer. We introduce a modified projection $\bm{\Pi}_k^{\rm m}$ for all $\bm v\in \bm H^{s}({\rm curl};\Omega)$ with $s>1/2$ and $\bm v_h\in \bm U_h$, such that:
\begin{align}\label{def1}
\bm{\Pi}_k^{\rm m}(\bm v,\bm v_h)=\bm{\mathcal B}_k\bm v+\nabla\sigma_h,
\end{align}
where
$\bm{\mathcal B}_k$ is some well-defined interpolation from $\bm H^s({\rm curl};\Omega)$ to $[\mathbb P_k(\mathcal{T}_h)]^3$, and \color{black} { $\sigma_h\in \mathbb P_{k}(\mathcal{T}_h)\cap H_0^1(\Omega)$  }\color{black} satisfies
\begin{align}\label{def2}
(\nabla\sigma_h,\nabla q_h)=(\bm{\Pi}_{h,k}^{\rm curl,c}(\bm v_h-\bm{\mathcal B}_k\bm v),\nabla q_h)\qquad  \color{black} { \forall\, q_h\in \mathbb P_{k}(\mathcal{T}_h)\cap H_0^1(\Omega).}\color{black}
\end{align}
With the above definition, we are ready to prove the following result.
\begin{lemma} For all $\bm v\in \bm H^{s}({\rm curl};\Omega)$ with $s>1/2$ and $\bm v_h\in \bm U_h$, there holds the orthogonality
	\begin{align}\label{oror}
	(\bm{\Pi}_{h,k}^{\rm curl,c}(\bm v_h-\bm{\Pi}_k^{\rm m}(\bm v,\bm v_h)),\nabla q_h)=0
	\qquad \color{black} { \forall\, q_h\in \mathbb P_{k}(\mathcal{T}_h)\cap H_0^1(\Omega),}\color{black}
	\end{align}
	and approximation properties
	\begin{align}
	\|\bm{\Pi}_k^{\rm m}(\bm v,\bm v_h)-\bm{\mathcal B}_k\bm v\|_0&\le  C \left(\|h_F^{1/2}\bm n\times[\![\bm v_h-\bm{\mathcal B}_k\bm v]\!]\|_{0,\mathcal{F}_h}+ \|\bm v_h-\bm{\mathcal B}_k\bm v\|_0 \right),\label{Bk}\\	
	\|\bm{\Pi}_{h,k}^{\rm curl,c}(\bm v_h-\bm{\Pi}_k^{\rm m}(\bm v,\bm v_h))\|_{0}&\le C\left(\|h_F^{1/2}\bm n\times[\![\bm v_h-\bm{\mathcal B}_k\bm v]\!]\|_{0,\mathcal{F}_h}+ \|\bm v_h-\bm{\mathcal B}_k\bm v\|_0 \right),\label{newpi1}\\
	\|\nabla\times((\bm{\Pi}_{h,k}^{\rm curl,c}\bm v_h-\bm{\Pi}_k^{\rm m}(\bm v,\bm v_h)))\|_{0}&\le C\left(\|h_F^{-1/2}\bm n\times[\![\bm v_h-\bm{\mathcal B}_k\bm v]\!]\|_{0,\mathcal{F}_h}  + \|\nabla_h\times(\bm v_h-\bm{\mathcal B}_k\bm v)\|_0 \right).\label{newpi2}
	\end{align}
\end{lemma}
\begin{proof}
	Using \eqref{def1}, \eqref{def2} and the fact that $\bm{\Pi}_{h,k}^{\rm curl,c}\nabla\sigma_h=\nabla\sigma_h$, it arrives at
	\begin{align*}
	(\bm{\Pi}_{h,k}^{\rm curl,c}(\bm v_h-\bm{\Pi}_k^{\rm m}(\bm v,\bm v_h)),\nabla q_h)&=	(\bm{\Pi}_{h,k}^{\rm curl,c}(\bm v_h-\bm{\mathcal B}_k\bm v)-\nabla \sigma_h,\nabla q_h)=0
	\\ \color{black} { \forall\, q_h\in \mathbb P_{k+1}(\mathcal{T}_h)\cap H_0^1(\Omega),}\color{black}
	\end{align*}
	which proves \eqref{oror}.
	We take $q_h=\sigma_h$ in \eqref{def2} to get
	\begin{align}
	\|\bm{\Pi}_k^{\rm m}(\bm v,\bm v_h)-\bm{\mathcal B}_k\bm v\|_0=
	\|\nabla\sigma_h\|_0\le \|\bm{\Pi}_{h,k}^{\rm curl,c}(\bm v_h-\bm{\mathcal B}_k\bm v)\|_0.\label{proof:bk}
	\end{align}
	Combining  \eqref{ic1}, \eqref{def1}, \eqref{proof:bk} and a triangle inequality gives
	\begin{align*}
	\|\bm{\Pi}_{h,k}^{\rm curl,c}(\bm v_h-\bm{\Pi}_k^{\rm m}(\bm v,\bm v_h))\|_{0}&\le 2\|\bm{\Pi}_{h,k}^{\rm curl,c}(\bm v_h-\bm{\mathcal B}_k\bm v)\|_0\nonumber\\
	&\le 2\left(\|\bm{\Pi}_{h,k}^{\rm curl,c}(\bm v_h-\bm{\mathcal B}_k\bm v)-(\bm v_h-\bm{\mathcal B}_k\bm v)\|_0
	+\|\bm v_h-\bm{\mathcal B}_k\bm v\|_0\right)\nonumber\\
	&\le C\left(\|h_F^{1/2}\bm n\times[\![\bm v_h-\bm{\mathcal B}_k\bm v]\!]\|_{0,\mathcal{F}_h}+\|\bm v_h-\bm{\mathcal B}_k\bm v\|_0\right),
	\end{align*}
	this implies \eqref{newpi1}. By using \eqref{proof:bk}, we can directly derive \eqref{Bk}. The proof of \eqref{newpi2} is similar to the proof of \eqref{newpi1}
\end{proof}

\color{black}

\subsubsection{H({\color{black}\rm div})-projection on domain surface}
For any $\bm{v}\in \bm{H}^s({\color{black}\rm curl};T)$ with $s>1/2$ and $T\in\mathcal{T}_h$,  we consider $\bm{\mathcal{P}}^{{\color{black}\rm curl}}_{k}\bm{v}$ restricting to the face $F$ such that
\begin{subequations}
	\begin{align}\label{314}
	\langle\!\langle\bm{\mathcal{P}}^{{\color{black}\rm curl}}_{k}\bm{v}\cdot\bm{t}_{FE},w_k\rangle\!\rangle_{E}
	=\langle\!\langle\bm{v}\cdot\bm{t}_{FE},w_k\rangle\!\rangle_{E}\qquad \forall\, w_k\in \mathbb{P}_k(E),\, E\subset\partial F,
	\end{align}
	where $\bm{t}_{FE}=\bm{n}_F\times\bm{n}_{FE}$, i.e., $\bm{t}_{FE}=\bm{t}_{{E}}$ or $\bm{t}_{FE}=-\bm{t}_{{E}}$. When $k\ge 2$, there holds
	\begin{align}\label{315}
	\langle \bm{\mathcal{P}}^{{\color{black}\rm curl}}_{k}\bm{v},\bm{w}_{k-1}  \rangle_{F}=\langle \bm{v},\bm{w}_{k-1}  \rangle_{F}\qquad\forall \,\bm{w}_{k-1}\in \bm{\mathcal{D}}_{k-1}(F).
	\end{align}
\end{subequations}
Since
\begin{align*}
\bm{v}\cdot\bm{t}_{FE}=\bm{v}\cdot(\bm{n}_{F}\times\bm{n}_{FE} )
= (\bm{v}\times\bm{n}_{F}) \cdot  \bm{n}_{E}
=-( \bm{n}_{F}\times\bm{v})\cdot  \bm{n}_{E} ,
\end{align*}
and
\begin{align*}
\bm{v}|_F&=(\bm{n}_F\times\bm{v})\times\bm{n}_F+(\bm{v}\cdot\bm{n}_F)\bm{n}_F,\qquad \bm{v}\cdot\bm{n}_F=0.
\end{align*}
Hence equations \eqref{314} and \eqref{315} can be rewritten as
\begin{align*}
\langle\!\langle\bm{\mathcal{P}}^{{\color{black}\rm curl}}_{k}( \bm{n}_{F}\times\bm{v})\cdot  \bm{n}_{E},w_k\rangle\!\rangle_{E}
=\langle\!\langle( \bm{n}_{F}\times\bm{v})\cdot \bm{n}_{E},w_k\rangle\!\rangle_{E}\qquad \forall\, w_k\in \mathbb{P}_k(E),\, E\subset\partial F,
\end{align*}
and when $k\ge 2$,
\begin{align*}
\langle \bm{\mathcal{P}}^{{\color{black}\rm curl}}_{k}(\bm{n}_F\times\bm{v}),\bm{n}_F\times\bm{w}_{k-1}  \rangle_{F}
=\langle (\bm{n}_F\times\bm{v}),\bm{n}_F\times\bm{w}_{k-1}  \rangle_{F}\quad\forall\, \bm{w}_{k-1}\in \bm{\mathcal{D}}_{k-1}(F).
\end{align*}
The operator $\bm{\mathcal{P}}^{{\color{black}\rm curl}}_{k}(\bm{n}_F\times\cdot)$ maps from space $\bm{H}^s({\color{black}\rm curl};T)$ to $\cup_{F\subset\partial T}[\mathbb{P}_k(F)]^3$ for each $T\in\mathcal{T}_h$. Actually, by denoting $\bm{\mathcal{P}}^{{\color{black}\rm div}}_{\Gamma,k}:=\bm{\mathcal{P}}^{{\color{black}\rm curl}}_{k}|_{\Gamma}$, we observe that $\bm{\mathcal{P}}^{{\color{black}\rm div}}_{\Gamma,k}(\bm{n}_F\times\bm{v})$ defines a $H({\color{black}\rm div})$-projection of $\bm{n}_F\times\bm{v}$ on domain surface $\Gamma$.

\subsection{HDG Method}

For any integers $k\ge 1$ and $m\in\{k-1,k\}$, we introduce  the following discrete spaces.
\begin{align*}
\bm{R}_{h}&=[\mathbb{P}_m(\mathcal{T}_h)]^3, \qquad
\bm{U}_{h}= [\mathbb{P}_k(\mathcal{T}_h)]^3,\qquad
\widehat{\bm{U}}_{h}=\{\widehat{\bm{v}}_h\in [\mathbb{P}_k(\mathcal{F}_h)]^2: \widehat{\bm{v}}_h\cdot\bm{n}|_{\mathcal{F}_h}=0\} ,\\
\widehat{\bm{U}}_{h}^{\widetilde{\bm{g}}}&=\{\widehat{\bm{v}}_h\in \widehat{\bm{U}}_{h}:
\bm{n} \times\widehat{\bm{v}}_h|_{\Gamma}=\bm{T}_{\Gamma,k}\widetilde{\bm{g}}  \} ,\qquad\widetilde{\bm{g}}=\bm{0},\bm{g}_T,\\
{P}_{h}&=\mathbb{P}_k(\mathcal{T}_h),\qquad
\widehat{P}_h=\mathbb{P}_k(\mathcal{F}_h),\qquad
\widehat{P}^0_{h}=\{\widehat{q}_h\in \widehat{P}_h:\widehat{q}_h|_{\Gamma}=0\},
\end{align*}
and $\bm{T}_{\Gamma,k}=\bm{\Pi}^{\partial}_k$ for the smooth case, $\bm{T}_{\Gamma,k}=\bm{\mathcal{P}}^{{\color{black}\rm div}}_{\Gamma,k}$ for the non-smooth case.

\begin{remark}
	Note that $\bm{\mathcal{P}}^{{\color{black}\rm div}}_{\Gamma,k}\widetilde{\bm{g}}$ with $\widetilde{\bm{g}}=\bm{g}_T$ in \eqref{source} is well-defined, since $\bm{g}_T=\bm{n} \times\bm{u}$ and Theorem \ref{reg} ensures $\bm{u}\in \bm{H}^{s}({\color{black}\rm curl};\Omega)$ with $s>1/2$. The calculation of $\bm{\mathcal{P}}^{{\color{black}\rm div}}_{\Gamma,k}\bm{g}_T$ can be done face by face (on $\Gamma$).
	We use $\bm{\mathcal{P}}^{{\color{black}\rm div}}_{\Gamma,k}\bm{g}_T$ to approximate $\bm{g}_T$ on surface $\Gamma$, and it could cancel some of the oscillation terms corresponding to $\bm g_T$ in the case where
	\begin{align}\label{eq:StabEst}
	\bm r\in [H^s(\Omega)]^3, \qquad\bm u\in  [H^{s}(\Omega)]^3,\qquad\bm g_T\in  \bm{H}^{-\frac{1}{2}}({\color{black}\rm div}_{\tau};\Gamma)\cap[H^{\delta}(\Gamma)]^3
	\end{align}
	with $s>1/2$ and $\delta\in(0,1]$. The numerical experiments in section~\ref{sec:NE} also demonstrate that $\bm{\mathcal{P}}^{{\color{black}\rm div}}_{\Gamma,k}\bm{g}_T$ is indeed a better choice for this case. However, in general, from Theorem \ref{reg} we see that only when $\delta>1$ in \eqref{eq:StabEst} could guarantee $\bm r\in [H^s(\Omega)]^3$ with $s>1/2$. With such regularity of $\bm g_T$, there is no need to use $\bm{\mathcal{P}}^{{\color{black}\rm div}}_{\Gamma,k}$ since no oscillation would be introduced by the $\bm g_T$ term on $\Gamma$.
\end{remark}
By introducing $\bm{r}=\nabla\times\bm{u}$ we can rewrite \eqref{source} as:  find $(\bm{r},\bm{u},p)$ that satisfies
\begin{eqnarray}\label{mixed}
\left\{
\begin{aligned}
\bm{r}-\nabla\times\bm{u}&=\bm{0}&\text{ in }\Omega,\\
\nabla\times\bm{r}+\nabla p&=\bm{f}&\text{ in }\Omega, \label{mix0source}\\
\nabla\cdot\bm{u}&=g&\text{ in }\Omega,  \\
\bm{n} \times \bm{u} &=\bm{g}_T&\text{ on }\Gamma,\\
p &=0&\text{ on } \Gamma.
\end{aligned}
\right.
\end{eqnarray}
The HDG finite element method for \eqref{mixed} in a compact form reads: for all $(\bm{s}_h,\bm{v}_h, \widehat{\bm{v}}_h,q_h, \widehat{q}_h) \in \bm{R}_h\times\bm{U}_h\times\widehat{\bm{U}}_h^{\bm{0}}\times {P}_h\times \widehat{{P}}^0_h$, find $(\bm{r}_h,\bm{u}_h,\widehat{\bm{u}}_h,p_h,\widehat{p}_h) \in \bm{R}_h\times\bm{U}_h\times\widehat{\bm{U}}_h^{\bm{g}_T}\times {P}_h\times \widehat{{P}}^0_h$ such that
\begin{subequations} \label{fem}
	\begin{align}
	a_h(\bm{r}_h,\bm{s}_h)+b_h(\bm{u}_h,\widehat{\bm{u}}_h;\bm{s}_h)&=0,\label{fhm1}\\
	b_h(\bm{v}_h,\widehat{\bm{v}}_h;\bm{r}_h)+c_h(p_h,\widehat{p}_h;\bm{v}_h)-s^u_h(  \bm{u}_h,\widehat{\bm{u}}_h; \bm{v}_h,\widehat{\bm{v}}_h )&=-(\bm{f},\bm{v}_h),\\
	c_h(q_h,\widehat{q}_h;\bm{u}_h)+s^p_h(p_h,\widehat{p}_h; q_h,\widehat{q}_h)&=(g,q_h),\label{fhm3}
	\end{align}
\end{subequations}
where
\begin{align*}
a_h(\bm{r}_h,\bm{s}_h)&=(\bm{r}_h,\bm{s}_h),\\
b_h(\bm{u}_h,\widehat{\bm{u}}_h;\bm{s}_h)&=-(\bm{u}_h,\nabla_h\times\bm{s}_h)-\langle \bm{n}\times\widehat{\bm{u}}_h,\bm{s}_h \rangle_{\partial\mathcal{T}_h},\\
c_h(q_h,\widehat q_h;\bm{u}_h)&=(\nabla_h\cdot\bm{u}_h, q_h)-\langle \bm{n}\cdot\bm{u}_h,\widehat{q}_h \rangle_{\partial\mathcal{T}_h},\\
s_h^u(  \bm{u}_h,\widehat{\bm{u}}_h; \bm{v}_h,\widehat{\bm{v}}_h )&=\langle h_F^{-1} \bm{n}\times(\bm{u}_h-\widehat{\bm{u}}_h), \bm{n}\times(\bm{v}_h-\widehat{\bm{v}}_h)\rangle_{\partial\mathcal{T}_h},\\
s^p_h(p_h,\widehat{p}_h; q_h,\widehat{q}_h)&=\langle h_F^{\alpha}(p_h-\widehat{p}_h),q_h-\widehat{q}_h \rangle_{\partial\mathcal{T}_h}\quad\text{ with }\alpha\in\{-1,1\}.
\end{align*}

To simplify the notation, we let
\begin{align*}
\bm{\sigma}_h &:=(\bm{r}_h ,\bm{u}_h ,\widehat{\bm{u}}_h ,p_h ,\widehat{p}_h ),\qquad\qquad\quad\;
\bm{\tau}_h :=(\bm{s}_h ,\bm{v}_h ,\widehat{\bm{v}}_h ,q_h ,\widehat{q_h }),\\
\bm{\Sigma}_h &:= \bm{R}_h\times\bm{U}_h\times\widehat{\bm{U}}_h\times P_h\times \widehat{P}_h,\qquad
\bm{\Sigma}_h^{\widetilde{\bm{g}}} := \bm{R}_h\times\bm{U}_h\times\widehat{\bm{U}}^{\widetilde{\bm{g}}}_h\times P_h\times \widehat{P}^0_h,\text{ with }
\widetilde{\bm{g}}=\bm{0},\bm{g}_T,
\end{align*}
and
\begin{align*}
B_h(\bm{\sigma}_h,\bm{\tau}_h)&:= a_h(\bm{r}_h,\bm{s}_h)+b_h(\bm{u}_h,\widehat{\bm{u}}_h;\bm{s}_h)+b_h(\bm{v}_h,\widehat{\bm{v}}_h;\bm{r}_h)\nonumber\\
&\quad+c_h(p_h,\widehat{p}_h;\bm{v}_h)-s^u_h(  \bm{u}_h,\widehat{\bm{u}}_h; \bm{v}_h,\widehat{\bm{v}}_h ) +c_h(q_h,\widehat{q}_h;\bm{u}_h)+s^p_h(p_h,\widehat{p}_h; q_h,\widehat{q}_h),\\
F_h(\bm{\tau}_h)&:=-(\bm{f},\bm{v}_h)+(g,q_h).
\end{align*}
The HDG method \eqref{fem} can then be rewritten as: find $\bm{\sigma}_h\in \bm{\Sigma}^{\bm{g}_T}_h$, such that
\begin{align}
B_h(\bm{\sigma}_h,\bm{\tau}_h)=F_h(\bm{\tau}_h)\qquad\forall\, \bm{\tau}_h\in\bm{\Sigma}^{\bm{0}}_h.\label{Bh_HDG}
\end{align}

\section{Stability Analysis}

We define semi-norms on the spaces $\bm{U}_h\times \widehat{\bm{U}}_h$ and $P_h\times \widehat{P}_h$ as follows:
\begin{align*}
\|(\bm{v},\widehat{\bm{v}})\|^2_{U} &:= \|(\bm{v},\widehat{\bm{v}})\|^2_{{\color{black}\rm curl}} +\|(\bm{v},\widehat{\bm{v}})\|^2_{{\color{black}\rm div}},\qquad
\|(q,\widehat{q})\|^2_{P} := \| h_T^{\frac{\alpha+1}{2}}\nabla_hq\|^2_0+ \|h_F^{\frac{\alpha}{2}}(q-\widehat{q})\|^2_{0,\partial\mathcal{T}_h},
\end{align*}
where
\begin{align*}
\|(\bm{v},\widehat{\bm{v}})\|^2_{{\color{black}\rm curl}}:=\|\nabla_{h}\times\bm{v}\|^2_0+\|h_F^{-1/2}(\bm{n}\times(\bm{v}-\widehat{\bm{v}}))\|^2_{0,\partial\mathcal{T}_h}, \\
\|(\bm{v},\widehat{\bm{v}})\|^2_{{\color{black}\rm div}}:=\|h_T^{\frac{1-\alpha}{2}}\nabla_{h}\cdot\bm{v}\|^2_0+\|h_F^{-\frac{\alpha}{2}}[\![\bm{n}\cdot\bm{v}]\!]\|^2_{0,\mathcal{F}^I_h}.
\end{align*}
Then we define the semi-norms on the space $\bm{\Sigma}_h$ as
\begin{align*}
\|\bm{\sigma} \|^2_{\bm{\Sigma}_h}&:=\|\bm{r} \|_0^2 +\|(\bm{u} ,\widehat{\bm{u} })\|^2_{U}+\|(p ,\widehat{p })\|^2_{P},\qquad
\|\bm{\tau} \|^2_{\bm{\Sigma}_h}:=\|\bm{s} \|_0^2+\|(\bm{v} ,\widehat{\bm{v} })\|^2_{U}+\|(q ,\widehat{q })\|^2_{P}.
\end{align*}

\begin{lemma}
	The semi-norm $\|(\cdot,\cdot)\|_{U}$ defines a norm on the space $\bm{U}_h\times\widehat{\bm{U}}^{\bm{0}}_h$.
\end{lemma}
\begin{proof} Let $(\bm{v}_h,\widehat{\bm{v}}_h)\in \bm{U}_h\times\widehat{\bm{U}}^{\bm{0}}_h$, it suffices to show that $\|(\bm{v}_h,\widehat{\bm{v}}_h)\|_{U}=0$ leads to $(\bm{v}_h,\widehat{\bm{v}}_h)=(\bm{0},\bm{0})$. It is obvious due to the fact that $\bm{n}\times\widehat{\bm{v}}_h=\bm{0}$ and $\bm{n}\cdot\widehat{\bm{v}}_h=0$ imply $\widehat{\bm{v}}_h=\bm{0}$, and $\nabla\cdot\bm v_h=0$, $\nabla\times\bm v_h=\bm 0$ imply $\bm v_h=\bm 0$ in $[H^1_0(\Omega)]^3$.
\end{proof}

\begin{lemma}
	The semi-norm $\|(\cdot,\cdot)\|_{P}$ defines a norm on ${P}_h\times\widehat{P}^0_h$, and for all $(q_h,\widehat{q}_h)\in {P}_h\times\widehat{{P}}^0_h$, there holds
	\begin{eqnarray}\label{relation}
	\|(q_h,\widehat{q}_h)\|^2_{P}\sim \| h_T^{\frac{\alpha+1}{2}}\nabla\mathcal{I}_k^cq_h\|^2_0 + \| h_F^{\frac{\alpha}{2}}(q_h-\widehat{q}_h)\|^2_{0,\partial\mathcal{T}_h}.
	\end{eqnarray}
\end{lemma}
\begin{proof}
	Utilizing the definition of $\|(\cdot,\cdot)\|_{P}$, a triangle inequality and the estimate \eqref{ick}, we get
	\begin{align*}
	\|(q_h,\widehat{q}_h)\|^2_{P}&=\|h_T^{\frac{\alpha+1}{2}}\nabla_h q_h\|^2_0+\|h_F^{\frac{\alpha}{2}}(q_h-\widehat{q}_h)\|^2_{0,\partial\mathcal{T}_h} \nonumber\\
	&\color{black}\le C \color{black}\left( \|h_T^{\frac{\alpha+1}{2}}(\nabla_h q_h-\nabla \mathcal{I}^c_kq_h)\|^2_0+\|h_T^{\frac{\alpha+1}{2}}\nabla \mathcal{I}^c_kq_h\|^2_0 +\|h_F^{\frac{\alpha}{2}}(q_h-\widehat{q}_h)\|^2_{0,\partial\mathcal{T}_h} \right)\nonumber\\
	&\color{black}\le C \color{black}\left( \|h_F^{\frac{\alpha}{2}}[\![q_h]\!]\|^2_{0,\mathcal{F}_h}+\|h_T^{\frac{\alpha+1}{2}}\nabla\mathcal{I}^c_kq_h\|^2_0+\|h_F^{\frac{\alpha}{2}}(q_h-\widehat{q}_h)\|^2_{0,\partial\mathcal{T}_h} \right)\nonumber\\
	&\color{black}\le C \color{black}\left( \|h_T^{\frac{\alpha+1}{2}}\nabla\mathcal{I}^c_kq_h\|^2_0+\|h_F^{\frac{\alpha}{2}}(q_h-\widehat{q}_h)\|^2_{0,\partial\mathcal{T}_h}\right) .
	\end{align*}
	Similarly, we have
	\begin{align*}
	\|h_T^{\frac{\alpha+1}{2}}\nabla\mathcal{I}^c_kq_h\|^2_0 &\color{black}\le C \color{black}\left( \|h_T^{\frac{\alpha+1}{2}}(\nabla\mathcal{I}^c_kq_h-\nabla_hq_h)\|^2_0+\|\nabla_hq_h\|^2_0\right)\nonumber\\
	&\color{black}\le C \color{black}\left( \|h_F^{\frac{\alpha}{2}}[\![q_h]\!]\|^2_{0,\mathcal{F}_h}+\|h_T^{\frac{\alpha+1}{2}}\nabla_hq_h\|^2_0 \right)\nonumber\\
	&\color{black}\le C \color{black}\left( \|h_F^{\frac{\alpha}{2}}(q_h-\widehat{q}_h)\|^2_{0,\partial\mathcal{T}_h}+\|\nabla_hq_h\|^2_0\right)\nonumber\\
	&= C\|(q_h,\widehat{q}_h)\|^2_{P}.
	\end{align*}
	Therefore, the estimate \eqref{relation} holds.
	
	Next, we prove that $\|(\cdot,\cdot)\|_P$ is a norm on ${P}_h\times\widehat{P}^0_h$.  For any $(q_h, \widehat{q}_h)\in {P}_h\times\widehat{P}^0_h$ such that $\|h_T^{\frac{\alpha+1}{2}}\nabla_hq_h\|^2_0
	+\|h_F^{\frac{\alpha}{2}}(q_h-\widehat{q}_h)\|^2_{0,\partial\mathcal{T}_h}=0$, we know $q_h$ is a piecewise constant, and $q_h=\widehat{q}_h$ on every face. Moreover, $q_h=\widehat{q}_h=0$ on boundary faces. Therefore, $q_h=\widehat{q}_h=0$. This completes the proof.
\end{proof}

\begin{theorem} [Discrete inf-sup condition] \label{Th45} The following stability results of $B_h$ hold
	\begin{subequations}
		\begin{align}
		\sup_{\bm 0\neq\bm{\tau}_h\in\bm{\Sigma}^{\bm{0}}_h}\frac{B_h(\bm{\sigma}_h,\bm{\tau}_h)}{\|\bm{\tau}_h\|_{\bm{\Sigma}_h}}
		\ge C   \|\bm{\sigma}_h\|_{\bm{\Sigma}_h},\label{lbb1}\\
		\sup_{\bm 0\neq\bm{\sigma}_h\in\bm{\Sigma}^{\bm{0}}_h}\frac{B_h(\bm{\sigma}_h,\bm{\tau}_h)}{\|\bm{\sigma}_h\|_{\bm{\Sigma}_h}}
		\ge C   \|\bm{\tau}_h\|_{\bm{\Sigma}_h}.
		\label{lbb2}
		\end{align}
	\end{subequations}
\end{theorem}

\begin{proof}
	Taking $\bm{\sigma}_h=(\bm{r}_h,\bm{u}_h,\widehat{\bm{u}}_h,p_h,\widehat{p}_h)\in\bm{\Sigma}^{\bm{0}}_h$, we establish \eqref{lbb1}--\eqref{lbb2} in five steps.

	\textbf{Step one:}
	
	Take
	$\bm{\tau}^1_h=(\bm{r}_h ,-\bm{u}_h ,-\widehat{\bm{u}}_h ,p_h ,\widehat{p}_h )\in\bm{\Sigma}^{\bm{0}}_h$, then by the definition of $\|\cdot\|_{\bm{\bm{\Sigma}}_h}$ and the definition of
	$B_h$ we have
	\begin{subequations}\label{B1}
		\begin{align}
		\|\bm{\tau}_h^1\|_{\bm{\Sigma}_h}=\|\bm{\sigma}_h\|_{\bm{\Sigma}_h},
		\end{align}
		and
		\begin{align}
		B_h(\bm{\sigma}_h,\bm{\tau}^1_h)=\|\bm{r}_h\|^2_0+\|h_F^{-1/2}\bm{n}\times(\bm{u}_h-\widehat{\bm{u}}_h)\|^2_{0,\partial\mathcal{T}_h}+ \|h_F^{\frac{\alpha}{2}}(p_h-\widehat{p}_h)\|^2_{0,\partial\mathcal{T}_h}.
		\end{align}
	\end{subequations}
	
	\textbf{Step two:}
	
	Taking $\bm{\tau}^2_h=(-\nabla_h\times\bm{u}_h,\bm{0},\bm{0},0,0) \in\bm{\Sigma}^{\bm{0}}_h$, then by the definition of $\|\cdot\|_{\bm{\bm{\Sigma}}_h}$ and the definition of
	$B_h$  we have
	\begin{subequations}\label{B2}
		\begin{align}
		\|\bm{\tau}_h^2\|_{\bm{\Sigma}_h}=\|\nabla_h\times\bm{u}_h\|_0\le\|\bm{\sigma}_h\|_{\bm{\Sigma}_h},
		\end{align}
		and using the definition of $B_h$, integration by parts, and an inverse inequality, we get
		\begin{align}
		B_h(\bm{\sigma}_h,\bm{\tau}^2_h)&=-(\bm{r}_h,\nabla_h\times\bm{u}_h) +\|\nabla_h\times\bm{u}_h\|^2_0-\langle\nabla_h\times\bm{u}_h,\bm{n}\times(\bm{u}_h- \widehat{\bm{u}}_h  )  \rangle_{\partial\mathcal{T}_h}\nonumber\\
		&\ge \frac{1}{2}\|\nabla_h\times\bm{u}_h\|^2_0-C_1\|\bm{r}_h\|^2_0-C_2 \|h_F^{-1/2}\bm{n}\times(\bm{u}_h-\widehat{\bm{u}}_h)\|^2_{0,\partial\mathcal{T}_h}.
		\end{align}
	\end{subequations}

	\textbf{Step three:}
	
	Let $r_h=h_T^{1-\alpha}\nabla_h\cdot\bm{u}_h$, $\widehat{r}_h=-h_F^{-\alpha}[\![\bm{n}\cdot \bm{u}_h]\!]$ on $\mathcal{F}_h^I$ and $\widehat{r}_h=0$ on $\Gamma$.
	We take $\bm{\tau}^3_h=(\bm{0},\bm{0},\bm{0},r_h,\widehat{r}_h)\in\bm{\Sigma}^{\bm{0}}_h$. Then by the definition of $\|\cdot\|_{\bm{\bm{\Sigma}}_h}$ and an inverse inequality we have
	\begin{subequations}\label{B3}
		\begin{align}
		\|\bm{\tau}^3_h\|_{\bm{\Sigma}_h}^2&=\|h_T^{\frac{\alpha+1}{2}}\nabla_hr_h\|^2_{0}+ \|h_F^{\frac{\alpha}{2}}(r_h-\widehat{r}_h)\|^2_{\partial\mathcal{T}_h}\nonumber\\
		&\color{black}\le C \color{black} \left( \| h_T^{\frac{1-\alpha}{2}}\nabla_h\cdot\bm{u}_h\|^2_0+ \|h_F^{-\frac{\alpha}{2}}[\![\bm{n}\cdot\bm{u}_h]\!]\|^2_{0,\mathcal{F}^I_h}\right)\nonumber\\
		&\le C\|\bm{\sigma}_h\|_{\bm{\Sigma}_h}^2,
		\end{align}
		Moreover, by the definition of $B_h$  we get
		\begin{align}
		B_h(\bm{\sigma}_h,\bm{\tau}^3_h)&=  \|h_T^{\frac{1-\alpha}{2}}\nabla\cdot\bm{u}_h\|^2_0+  \| h_F^{-\frac{\alpha}{2}}[\![\bm{n}\cdot\bm{u}_h]\!]\|^2_{0,\mathcal{F}_h^I} + \langle h_F^{\alpha} (p_h-\widehat{p}_h), r_h-\widehat{r}_h\rangle_{\partial\mathcal{T}_h} \nonumber\\
		&\ge\frac{1}{2}\|(\bm{u}_h,\widehat{\bm{u}}_h)\|^2_{{\color{black}\rm div}}-C_3 \|h_F^{\frac{\alpha}{2}}(p_h-\widehat{p}_h )\|_{0,\mathcal{F}_h}.
		\end{align}
	\end{subequations}

	\textbf{Step four:}
	
	We then take $\bm{\tau}^4_h=(\bm{0},-h_T^{\alpha+1}\nabla\mathcal{I}^c_k p_h,-h_F^{\alpha+1}\bm{n}\times\nabla\mathcal{I}^c_k p_h\times\bm{n},0,0)\in\bm{\Sigma}^{\bm{0}}_h$, and consider the cases where $\alpha=-1$ and $\alpha=1$ separately.
	
	{\bf (i)} When $\alpha=-1$, by the definition of $\|\cdot\|_{\bm{\bm{\Sigma}}_h}$ and an inverse inequality, we have
	%\
	\begin{subequations}\label{B41}
		\begin{align}
		\|\bm{\tau}^4_h\|_{\bm{\Sigma}_h}^2
		&=\|h_T^{\frac{1-\alpha}{2}}\nabla_h\cdot h_T^{\alpha+1}\nabla\mathcal{I}^c_k p_h\|^2_0+ \|h_F^{-\frac{\alpha}{2}}[\![h_T^{\alpha+1}\nabla\mathcal{I}^c_k p_h\cdot\bm{n}]\!]\|_{0,\mathcal{F}^I_h}^2\le C \color{black}  \|h_T^{\frac{\alpha+1}{2}}\nabla\mathcal{I}^c_k p_h\|^2_0\le C \color{black} \|\bm{\sigma}_h\|_{\bm{\Sigma}_h}^2.
		\end{align}
		By the definition of $B_h$ we have
		\begin{align}
		B_h(\bm{\sigma}_h,\bm{\tau}^4_h)&=(h_T^{\alpha+1}\nabla\mathcal{I}^c_k p_h,\nabla_h\times\bm{r}_h)+\langle \bm{n}\times h_T^{\alpha+1}\nabla\mathcal{I}^c_k p_h,\bm{r}_h\rangle_{\partial\mathcal{T}_h}\nonumber\\
		&\quad
		-  (\nabla_h\cdot h_T^{\alpha+1}\nabla\mathcal{I}^c_k p_h, p_h)+ \langle \bm{n}\cdot h_T^{\alpha+1}\nabla\mathcal{I}^c_k p_h,\widehat{p}_h \rangle_{\partial\mathcal{T}_h}\nonumber\\
		&=(h_T^{\alpha+1}\nabla_h p_h, \nabla_h p_h)+ (h_T^{\alpha+1}(\nabla\mathcal{I}^c_k p_h- \nabla_hp_h),\nabla_h p_h)
		+ \langle  h_T^{\alpha+1}\bm{n}\cdot\nabla\mathcal{I}^c_k p_h,\widehat{p}_h -p_h \rangle_{\partial\mathcal{T}_h}\nonumber\\
		&\ge \frac{1}{2}\|h_T^{\frac{\alpha+1}{2}}\nabla_h p_h\|^2_0-C_4  \| h_F^{\frac{\alpha}{2}}(p_h-\widehat{p}_h)\|^2_{0,\partial\mathcal{T}_h}.
		\end{align}
	\end{subequations}

	{\bf (ii)} When $\alpha=1$, by the definition of $\|\cdot\|_{\bm{\bm{\Sigma}}_h}$ and an inverse inequality, we have
	\begin{subequations}\label{B42}
		\begin{align}
		\|\bm{\tau}^4_h\|_{\bm{\Sigma}_h}^2&= \|h_T^{\frac{1-\alpha}{2}}\nabla_h\cdot h_T^{\alpha+1}\nabla\mathcal{I}^c_k p_h\|^2_0+ \|h_F^{-\frac{\alpha}{2}}[\![h_T^{\alpha+1}\nabla\mathcal{I}^c_k p_h\cdot\bm{n}]\!]\|_{0,\mathcal{F}^I_h}^2\nonumber\\
		&\qquad+\|h_F^{-1/2}(h_T^{\alpha+1}\bm{n}\times\nabla\mathcal{I}^c_k p_h-h_F^{\alpha+1}\bm{n}\times\nabla\mathcal{I}^c_k p_h)\|^2_{0,\partial\mathcal{T}_h}\nonumber\\
		&\color{black}\le C \color{black}  \|h_T^{\frac{\alpha+1}{2}}\nabla\mathcal{I}^c_k p_h\|^2_0\nonumber\\
		&\color{black}\le C \color{black} \|\bm{\sigma}_h\|_{\bm{\Sigma}_h}^2.
		\end{align}
		From the definition of $B_h$, we obtain
		\begin{align}\label{eq:Bh6Est}
		B_h(\bm{\sigma}_h,\bm{\tau}^4_h)&=(h_T^{\alpha+1}\nabla\mathcal{I}^c_k p_h,\nabla_h\times\bm{r}_h)+\langle \bm{n}\times h_F^{\alpha+1}\nabla\mathcal{I}^c_k p_h,\bm{r}_h\rangle_{\partial\mathcal{T}_h}
		\nonumber\\
		&\quad-  (\nabla_h\cdot h_T^{\alpha+1}\nabla\mathcal{I}^c_k p_h, p_h)+ \langle \bm{n}\cdot h_T^{\alpha+1}\nabla\mathcal{I}^c_k p_h,\widehat{p}_h \rangle_{\partial\mathcal{T}_h}\nonumber\\
		&\quad+\langle h_F^{-1}(h_T^{\alpha+1}\bm{n}\times\nabla\mathcal{I}^c_k p_h-h_F^{\alpha+1}\bm{n}\times\nabla\mathcal{I}^c_k p_h),
		\bm{n}\times(\bm{v}_h-\widehat{\bm{v}}_h) \rangle_{\partial\mathcal{T}_h}\nonumber\\
		&=
		(\nabla_h\times (h_T^{\alpha+1}\nabla\mathcal{I}^c_k p_h),\bm{r}_h)
		+
		(h_T^{\alpha+1}\nabla_h p_h, \nabla_h p_h)+ (h_T^{\alpha+1}(\nabla\mathcal{I}^c_k p_h- \nabla_hp_h),\nabla_h p_h)
		\nonumber\\
		&\quad+  \langle  h_T^{\alpha+1}\bm{n}\cdot\nabla\mathcal{I}^c_k p_h,\widehat{p}_h -p_h \rangle_{\partial\mathcal{T}_h}\nonumber\\
		&\quad+\langle h_F^{-1}(h_T^{\alpha+1}\bm{n}\times\nabla\mathcal{I}^c_k p_h-h_F^{\alpha+1}\bm{n}\times\nabla\mathcal{I}^c_k p_h), \bm{n}\times(\bm{v}_h-\widehat{\bm{v}}_h) \rangle_{\partial\mathcal{T}_h}\nonumber\\
		&\ge \frac{1}{2}\|h_T^{\frac{\alpha+1}{2}}\nabla_h p_h\|^2_0-C_4  \| h_F^{\frac{\alpha}{2}}(p_h-\widehat{p}_h)\|^2_{0,\partial\mathcal{T}_h}-C_5 \| h_F^{-1/2}\bm{n}\times(\bm{v}_h-\widehat{\bm{v}}_h) \|^2_{0,\partial\mathcal{T}_h}-C_6\|\bm{r}_h\|^2_0.
		\end{align}
	\end{subequations}

	\textbf{Step five:}
	
	Take $C_0=\max(C_1+C_6,C_2+C_5,C_3+C_4)+1$ and
	$\bm{\tau}_h=C_0\bm{\tau}_1+\bm{\tau}_2+\bm{\tau}_3+\bm{\tau}_4$.
	Then combining \eqref{B1}--\eqref{B42}, we arrive at
	\begin{eqnarray*}
		\|\bm{\tau}_h\|_{\bm{\Sigma}_h}\color{black}\le C \color{black} \|\bm{\sigma}_h\|_{\bm{\Sigma}_h},
	\end{eqnarray*}
	and
	\begin{align*}
	B_h(\bm{\sigma}_h,\bm{\tau}_h) \ge C \left(\|\bm{r}_h\|_0^2+\|(\bm{u}_h,\widehat{\bm{u}}_h)\|^2_{U}+\|(p_h,\widehat{p}_h)\|^2_{P}\right)=C\|\bm{\sigma}_h\|^2_{\bm{\Sigma}_h}.
	\end{align*}
	Finally, the last two inequalities lead to
	\begin{align*}
	B_h(\bm{\sigma}_h,\bm{\tau}_h)\ge C  \|\bm{\sigma}_h\|_{\bm{\Sigma}_h}\|\bm{\tau}_h\|_{\bm{\Sigma}_h},
	\end{align*}
	which implies \eqref{lbb1}. Since $B_h$ is symmetric, \eqref{lbb2} also holds.
\end{proof}

The next conclusion is a direct consequence of Theorem~\ref{Th45}.
\begin{corollary}
	The HDG method \eqref{Bh_HDG} admits a unique solution $\bm{\sigma}_h\in\bm{\Sigma}_h^{\bm{0}}$.
\end{corollary}

\section{Error estimates}

\subsection{High Regularity Case}

In this subsection, we assume that $\bm{r},\bm{u}$ and $p$ are smooth functions (i.e., the solutions have high regularities). Let $\bm{T}_{\Gamma,k}=\bm{\Pi}^{\partial}_k$, $\alpha\in \{-1,1\}$, and
\begin{align*}
t:=s-\frac{1+\alpha}{2}\in \left(\frac{1}{2},k+1\right], \text{ where} \,s\in \left(\frac{3}{2},k+1\right].
\end{align*}

\subsubsection{Primary error estimates}
\begin{lemma} \label{Jh}Let $(\bm{r},\bm{u},p)$ be the solution of \eqref{mixed}, $\bm{\sigma}$ and $\bm{\mathcal{J}}_h\bm{\sigma}$ be defined as
	\begin{align*}
	\bm{\sigma}&:=(\bm{r},\bm{u},\widehat{\bm{u}},p,\widehat{p}), \,\,\,\text{ where } (\widehat{\bm{u}},\widehat{p})=(\bm{u},p) \,\text{on}\, \mathcal{F}_h,\qquad
	\bm{\mathcal{J}}_h\bm{\sigma}:=(\bm{\Pi}^{o}_m\bm{r}, \bm{\Pi}^{o}_k\bm{u}, \bm n\times\bm{\Pi}^{\partial}_k\bm{u}\times\bm n, \Pi^{o}_kp, \Pi^{\partial}_kp).
	\end{align*}
	Then we have the equation
	\begin{align}\label{510}
	B_h(\bm{\mathcal{J}}_h\bm{\sigma},\bm{\tau}_h)=F_h(\bm{\tau}_h)+E^{\mathcal{J}}_h(\bm{\sigma};\bm{\tau}_h)\quad\forall\, \bm{\tau}_h\in\bm{\Sigma}^{\bm{0}}_h,
	\end{align}
	where
	\begin{align}\label{eq:EJh}
	E^{\mathcal{J}}_h(\bm{\sigma};\bm{\tau}_h)&= -\langle \bm{n}\times(\widehat{\bm{v}}_h-\bm{v}_h),\bm{\Pi}^{o}_m\bm{r}-\bm{r} \rangle_{\partial\mathcal{T}_h} -\langle h_F^{-1} \bm{n}\times(  \bm{\Pi}^{o}_k\bm{u}-\bm{u}        ), \bm{n}\times(\bm{v}_h-\widehat{\bm{v}}_h  )  \rangle_{\partial\mathcal{T}_h}\nonumber\\
	&\quad-\langle \bm{n}\cdot(\bm{\Pi}^{o}_k\bm u-\bm{u}),\widehat{q}_h -q_h \rangle_{\partial\mathcal{T}_h}+\langle h_F^{\alpha}(\Pi^{o}_kp-p),q_h-\widehat{q}_h \rangle_{\partial\mathcal{T}_h}.
	\end{align}
\end{lemma}

\begin{proof} Utilizing the definitions of $a_h$ and $b_h$ we have
	\begin{eqnarray}\label{eq:abh}
	&&a_h(\bm{\Pi}^{o}_m\bm{r},\bm{s}_h)+b_h(\bm{\Pi}^{o}_k\bm{u},\bm n\times\bm{\Pi}^{\partial}_k\bm{u}\times\bm n;\bm{s}_h)\nonumber\\
	&&\quad\quad= (\bm{\Pi}^{o}_m\bm{r},\bm{s}_h)-(\bm{\Pi}^{o}_k\bm{u},\nabla_h\times\bm{s}_h)-\langle \bm{n}\times\bm{\Pi}^{\partial}_k\bm{u},\bm{s}_h \rangle_{\partial\mathcal{T}_h} \nonumber\\
	&&\quad\quad= (\bm{r},\bm{s}_h)-(\nabla\times\bm{u},\bm{s}_h) \nonumber\\
	&&\quad\quad= 0.
	\end{eqnarray}
	Using the definitions of $b_h$, $c_h$ and $s_h^u$ one can get
	\begin{eqnarray}\label{eq:bcsh}
	&&b_h(\bm{v}_h,\widehat{\bm{v}}_h;\bm{\Pi}^{o}_m\bm{r})+c_h(\Pi^{o}_kp,\Pi^{\partial}_kp;\bm{v}_h) -s^u_h(  \bm{\Pi}^{o}_k\bm{u},\bm n\times\bm{\Pi}^{\partial}_k\bm{u}\times\bm n; \bm{v}_h,\widehat{\bm{v}}_h ) \nonumber\\
	&&\quad\quad=-(\bm{v}_h,\nabla_h\times\bm{\Pi}^{o}_m\bm{r})-\langle \bm{n}\times\widehat{\bm{v}}_h,\bm{\Pi}^{o}_m\bm{r} \rangle_{\partial\mathcal{T}_h} +(\nabla_h\cdot\bm{v}_h, \Pi^{o}_kp) \nonumber\\
	&&\quad\qquad-\langle \bm{n}\cdot\bm{v}_h,\Pi^{\partial}_kp \rangle_{\partial\mathcal{T}_h}-\langle h_F^{-1} \bm{n}\times(  \bm{\Pi}^{o}_k\bm{u}-\bm{u}        ), \bm{n}\times(\bm{v}_h-\widehat{\bm{v}}_h  )  \rangle_{\partial\mathcal{T}_h}\nonumber\\
	&&\quad\quad=-(\bm{v}_h,\nabla\times\bm{r})-\langle \bm{n}\times(\widehat{\bm{v}}_h-\bm{v}_h),\bm{\Pi}^{o}_m\bm{r}-\bm{r} \rangle_{\partial\mathcal{T}_h}-(\bm{v}_h, \nabla p)
	\nonumber\\
	&&\qquad\quad-\langle h_F^{-1} \bm{n}\times(  \bm{\Pi}^{o}_k\bm{u}-\bm{u}        ), \bm{n}\times(\bm{v}_h-\widehat{\bm{v}}_h  )  \rangle_{\partial\mathcal{T}_h}
	\nonumber\\
	&&\quad\quad=-(\bm{f},\bm{v}_h)-\langle \bm{n}\times(\widehat{\bm{v}}_h-\bm{v}_h),\bm{\Pi}^{o}_m\bm{r}-\bm{r} \rangle_{\partial\mathcal{T}_h}-\langle h_F^{-1} \bm{n}\times(  \bm{\Pi}^{o}_k\bm{u}-\bm{u}        ), \bm{n}\times(\bm{v}_h-\widehat{\bm{v}}_h  )  \rangle_{\partial\mathcal{T}_h}.
	\end{eqnarray}
	From the definitions of $c_h$ and $s_h^p$ one can get
	\begin{align}\label{csh}
	&c_h(q_h,\widehat{q}_h;\bm{\Pi}^{o}_k\bm{u})+s^p_h(\Pi^{o}_kp,\Pi^{\partial}_kp; q_h,\widehat{q}_h) \nonumber\\
	&\qquad=(\nabla_h\cdot\bm{\Pi}^{o}_k\bm{u}, q_h)-\langle \bm{n}\cdot\bm{\Pi}^{o}_k\bm{u},\widehat{q}_h \rangle_{\partial\mathcal{T}_h}+\langle h_F^{\alpha}(\Pi^{o}_kp-\Pi^{\partial}_kp),q_h-\widehat{q}_h \rangle_{\partial\mathcal{T}_h} \nonumber\\
	&\qquad=(\nabla\cdot\bm{u}, q_h)-\langle \bm{n}\cdot(\bm{\Pi}^{o}_k\bm u-\bm{u}),\widehat{q}_h -q_h \rangle_{\partial\mathcal{T}_h}+\langle h_F^{\alpha}(\Pi^{o}_kp-p),q_h-\widehat{q}_h \rangle_{\partial\mathcal{T}_h} \nonumber\\
	&\qquad=(g,q_h)-\langle \bm{n}\cdot(\bm{\Pi}^{o}_k\bm u-\bm{u}),\widehat{q}_h -q_h \rangle_{\partial\mathcal{T}_h}+\langle h_F^{\alpha}(\Pi^{o}_kp-p),q_h-\widehat{q}_h \rangle_{\partial\mathcal{T}_h}.
	\end{align}
	Therefore, equality \eqref{510} directly follows from the definition of $B_h$, and \eqref{eq:abh}--\eqref{csh}.
\end{proof}

\begin{lemma} \label{lemma53}
	Let $\bm{\sigma}$ be as defined as in Lemma \ref{Jh} and $(\bm{r},\bm{u},p)\in [H^{s-1}(\Omega)]^3\times[H^{s}(\Omega)]^3\times H^{t}(\Omega)$ with $s\in (3/2,k+1]$ and $t=s-(\alpha+1)/2$, then there holds
	\begin{equation}\label{eq:EJh_Est}
	E^{\mathcal{J}}_h(\bm{\sigma};\bm{\tau}_h)\color{black}\le C \color{black} h^{s-1}(\|\bm{r}\|_{s-1}+\|\bm{u}\|_{s}+\|p\|_{t} )\|\bm{\tau}_h\|_{\bm{\Sigma}_h}.
	\end{equation}
\end{lemma}
\begin{proof}
	To simplify the notation, we define
	\begin{align}
	E^{\mathcal{J}}_1&= -\langle \bm{n}\times(\widehat{\bm{v}}_h-\bm{v}_h),\bm{\Pi}^{o}_m\bm{r}-\bm{r} \rangle_{\partial\mathcal{T}_h},\qquad
	E^{\mathcal{J}}_2=-\langle h_F^{-1} \bm{n}\times(  \bm{\Pi}^{o}_k\bm{u}-\bm{u}        ), \bm{n}\times(\bm{v}_h-\widehat{\bm{v}}_h  )  \rangle_{\partial\mathcal{T}_h},\nonumber\\
	E^{\mathcal{J}}_3&=-\langle \bm{n}\cdot(\bm{\Pi}^{o}_k\bm u-\bm{u}),\widehat{q}_h -q_h \rangle_{\partial\mathcal{T}_h},\qquad
	E^{\mathcal{J}}_4=\langle h_F^{\alpha}(\Pi^{o}_kp-p),q_h-\widehat{q}_h \rangle_{\partial\mathcal{T}_h}.\nonumber
	\end{align}
	We then estimate $E^{\mathcal{J}}_i$ term by term as follows.
	\begin{align}
	|E^{\mathcal{J}}_1|&\color{black}\le C \color{black} \sum_{T\in\mathcal{T}_h}\|\bm{n}\times(\widehat{\bm{v}}_h-\bm{v}_h)\|_{0,\partial T}h_T^{s-3/2}\|\bm{r}\|_{s-1,T}\color{black}\le C \color{black} h^{s-1}\|\bm{r}\|_{s-1}\|\bm{\tau}_h\|_{\bm{\Sigma}_h},\label{eq:E1h_est}\\
	|E^{\mathcal{J}}_2|&\color{black}\le C \color{black} \sum_{T\in\mathcal{T}_h}\|\bm{n}\times(\widehat{\bm{v}}_h-\bm{v}_h)\|_{0,\partial T}h_T^{s-3/2}\|\bm{u}\|_{s,T}\color{black}\le C \color{black} h^{s-1}\|\bm{u}\|_{s}\|\bm{\tau}_h\|_{\bm{\Sigma}_h},\label{eq:E2h_est}\\
	|E^{\mathcal{J}}_3|&\color{black}\le C \color{black} \sum_{T\in\mathcal{T}_h}\|q_h-\widehat{q}_h\|_{0,\partial T}h_T^{s-1/2}\|\bm{u}\|_{s,T}\color{black}\le C \color{black} h^{s-1}\|\bm{u}\|_{s}\|\bm{\tau}_h\|_{\bm{\Sigma}_h},\label{eq:E3h_est}\\
	|E^{\mathcal{J}}_4|&\color{black}\le C \color{black} \sum_{T\in\mathcal{T}_h}\|q_h-\widehat{q}_h\|_{0,\partial T}h_T^{t-1/2+\alpha}\|p\|_{t,T}\color{black}\le C \color{black} h^{t+(\alpha-1)/2}\|p\|_{t}\|\bm{\tau}_h\|_{\bm{\Sigma}_h}=h^{s-1}\|p\|_{t}\|\bm{\tau}_h\|_{\bm{\Sigma}_h}.\label{eq:E4h_est}
	\end{align}
	Hence, the estimation \eqref{eq:EJh_Est} directly follows from \eqref{eq:EJh}, \eqref{eq:E1h_est}--\eqref{eq:E4h_est}.
\end{proof}

\begin{theorem} \label{l53}  Let $\bm{\sigma}$ be defined as in Lemma \ref{Jh}
	and $(\bm{r},\bm{u},p)\in [H^{s-1}(\Omega)]^3\times[H^{s}(\Omega)]^3\times H^{t}(\Omega)$ with $s\in (3/2,k+1]$ and $t=s-(\alpha+1)/2>1/2$, let $\bm\sigma_h$ be the solution of \eqref{Bh_HDG}, then there holds
	\begin{align*}
	\|\bm{\sigma}_h-\bm{\sigma}\|_{\bm{\Sigma}_h}\color{black}\le C \color{black} h^{s-1}(\|\bm{r}\|_{s-1}+\|\bm{u}\|_{s}+\|p\|_{t} ).
	\end{align*}
\end{theorem}

\begin{proof}
	From Theorem \ref{Th45} and, \eqref{Bh_HDG}, Lemma \ref{Jh}, and Lemma \ref{lemma53}, we have
	\begin{align*}
	\|\bm{\sigma}_h-\bm{\mathcal{J}}_h\bm{\sigma}\|_{\bm{\Sigma}_h} &\color{black}\le C \color{black} \sup_{\bm 0\neq\bm{\tau}_h\in\bm{\Sigma}_h^{\bm 0}}\frac{B_h( \bm{\sigma}_h-{\bm{\mathcal{J}}_h}\bm{\sigma},\bm{\tau}_h  )}{\|\bm{\tau}_h\|_{\bm{\Sigma}_h}}=C \sup_{\bm 0\neq\bm{\tau}_h\in\bm{\Sigma}_h^{\bm 0}}\frac{E^{\mathcal{J}}_h( \bm{\sigma}_h,\bm{\tau}_h  )}{\|\bm{\tau}_h\|_{\bm{\Sigma}_h}}\nonumber\\ 
	&\le C \color{black} h^{s-1}(\|\bm{r}\|_{s-1}+\|\bm{u}\|_{s}+\|p\|_{t} ),
	\end{align*}
	which combines with a triangle inequality to get the result.
\end{proof}

\subsubsection{$L^2$ error estimates}
Assume $\bm\Theta\in \bm{H}({\color{black}\rm div} ^0;\Omega)$. We first introduce the dual problem:
\begin{eqnarray}
\left\{
\begin{aligned}\label{eq:dualProb}
\bm{r}^d-  \nabla\times \bm{u}^d&=0,&\text{ in }\Omega,\\
\nabla\times \bm{r}^d+\nabla p^d&=\color{black}\bm {\Theta}, &\text{ in }\Omega, \label{duall}\\
\nabla\cdot\bm{u}^d&=\color{black}\Lambda,  &\text{ in }\Omega,  \\
\bm{n} \times \bm{u}^d &=\bm{0}, &\text{ on }\Gamma,\\
p^d&=0, &\text{ on }\Gamma.
\end{aligned}
\right.
\end{eqnarray}
Moreover, we assume that the solution $(\bm{r}^d, \bm{u}^d)$ of the dual problem \eqref{eq:dualProb} satisfies
\begin{align} \color{black} \label{reg1}
\|\bm{r}^d\|_{\beta}+\|\bm{u}^d\|_{1+\beta}\le C\left( \color{black} \|\color{black}\bm {\Theta}\|_0+\|\Lambda\|_0\right).
\end{align}
for some fixed $\beta\in \left(1/2,1\right]$. Notice that $p^d=0$ due to $\bm\Theta\in \bm{H}({\color{black}\rm div} ^0;\Omega)$.

\begin{lemma}\label{lemma55} Let $(\bm r,\bm u, p)$ and $(\bm r^d,\bm u^d,p^d)$ be the solution of \eqref{mix0source} and \eqref{duall}, respectively, let $\bm\sigma_h$ be the solution of \eqref{Bh_HDG},  then there holds
	\begin{align}
	|E^{\mathcal{J}}_h(\bm{\sigma};\bm{\mathcal{J}}_h\bm{\sigma}^d)|&\color{black}\le C \color{black} h^{s-1+\beta}(\|\bm{r}\|_{s-1}+\|\bm{u}\|_{s} ) \bm \left(\|\bm{\Theta}\|_0+\|\Lambda\|_0\right),\label{eq:EJSJh} \\
	|E^{\mathcal{J}}_h(\bm{\sigma}^d; {\bm{\mathcal{J}}_h}\bm{\sigma}-\bm{\sigma}_h)|&\color{black}\le C \color{black} h^{s-1+\beta}(\|\bm{r}\|_{s-1}+\|\bm{u}\|_{s}+\|p\|_{t} )\left(\|\bm{\Theta}\|_0+\|\Lambda\|_0\right).\label{eq:EJSJhD}
	\end{align}
	Here
	$\bm{\sigma}^d:=(\bm{r}^d,\bm{u}^d,\widehat{\bm{u}}^d,p^d,\widehat{p}^d)$, with $(\widehat{\bm{u}}^d,\widehat{p}^d)=(\bm{u}^d,p^d)$ on $\mathcal{F}_h$.
\end{lemma}
\begin{proof} By the definition of $E^{\mathcal{J}}_h$ (cf. \eqref{eq:EJh}), we have
	\begin{align*}
	E^{\mathcal{J}}_h(\bm{\sigma};\bm{\mathcal{J}}_h\bm{\sigma}^d)&\le |\langle \bm{n}\times(\bm{\Pi}^{\partial}_k\bm{u}^d-\bm{\Pi}^{o}_k\bm{u}^d),\bm{\Pi}^{o}_m\bm{r}-\bm{r} \rangle_{\partial\mathcal{T}_h} | \nonumber\\
	&\quad+|\langle h_F^{-1} \bm{n}\times(  \bm{\Pi}^{o}_k\bm{u}-\bm{u} ), \bm{n}\times(\bm{\Pi}^{o}_k\bm{u}^d-\bm{\Pi}^{\partial}_k\bm{u}^d  )  \rangle_{\partial\mathcal{T}_h}|\nonumber\\
	&\color{black}\le C \color{black} h^{s-1+\beta}(\|\bm{r}\|_{s-1}+\|\bm{u}\|_{s} )\|\bm{u}^d\|_{1+\beta} \nonumber\\
	&\color{black}\le C \color{black} h^{s-1+\beta}(\|\bm{r}\|_{s-1}+\|\bm{u}\|_{s} ) \left(\|\bm{\Theta}\|_0+\|\Lambda\|_0\right).
	\end{align*}
	Moreover, it follows from Lemma \ref{lemma53} that
	\begin{align*}
	|E^{\mathcal{J}}_h(\bm{\sigma}^d; {\bm{\mathcal{J}}_h}\bm{\sigma}-\bm{\sigma}_h)|&\color{black}\le  C \color{black} h^{\beta}(\|\bm{r}^d\|_{\beta}+\|\bm{u}^d\|_{\beta+1} )
	\| {\bm{\mathcal{J}}_h}\bm{\sigma}-\bm{\sigma}_h\|_{\bm{\Sigma}_h}\nonumber\\
	&\color{black}\le C \color{black} h^{s-1+\beta}
	\left(\|\bm{\Theta}\|_0+\|\Lambda\|_0\right)(\|\bm{r}\|_{s-1}+\|\bm{u}\|_{s}+\|p\|_{t} ).
	\end{align*}
	This completes the proof.
\end{proof}

\begin{theorem}\label{L^2}
	Let $(\bm{r},\bm{u},p)\in [H^{s-1}(\Omega)]^3\times[H^{s}(\Omega)]^3\times H^{t}(\Omega)$ with $s\in \left(3/2,k+1\right]$, $t=s-(\alpha+1)/2>1/2$ be the solution of \eqref{mixed},  let $(\bm r_h,\bm u_h,\widehat{\bm u}_h,p_h,\widehat p_h)$ be the solution of \eqref{fem}, then there holds
	\color{black}
	\begin{align}\label{eq:PkUP}
	\|\bm{u}-\bm{u}_h\|_0\color{black}+\|\Pi_k^op-p_h\|_0&\le C  \left(h^{s-1+(\sigma+\beta)/2}
	+h^{s-1+\beta}
	\right)
	(\|\bm{r}\|_{s-1}+\|\bm{u}\|_{s}+\|p\|_{t} ),
	\end{align}
	where $\sigma\in (1/2,1]$ is defined in Lemma \ref{embed}, and $\beta\in (1/2,1]$ is defined in \eqref{reg1}.
\end{theorem}
\begin{proof}
	\color{black}
	Take $\Lambda=\Pi_k^op-p_h$ in \eqref{duall} and let $\bm\Theta\in \bm H({\rm curl};\Omega)\cap \bm H({\rm div};\Omega)$ be the solution of
	\begin{subequations}
		\begin{align*}
		\nabla\times\bm\Theta&=\nabla\times ( \bm{\Pi}_{h,k}^{\rm curl,c}(\bm u_h-\bm{\Pi}_k^{\rm m}(\bm u,\bm u_h)))&\text{in }\Omega,\\
		\nabla\cdot\bm\Theta&=0&\text{in }\Omega,\\
		\bm n\times\bm\Theta&=\bm 0&\text{on }\Gamma.
		\end{align*}
	\end{subequations}
	Due to \eqref{oror} and the result in \cite[Lemma 4.5]{MR2009375} one has
	\begin{align}
	%\|\bm\Theta\|_{\sigma}&\le C\|\nabla\times ( \bm{\Pi}_{h,k}^{\rm curl,c}(\bm u_h-\bm{\Pi}_k^{\rm m}(\bm u,\bm u_h) ))\|_0,\label{co-01}\\
	\|\bm\Theta-( \bm{\Pi}_{h,k}^{\rm curl,c}(\bm u_h-\bm{\Pi}_k^{\rm m}(\bm u,\bm u_h) ))\|_{0}&\le Ch^{\sigma}\|\nabla\times ( \bm{\Pi}_{h,k}^{\rm curl,c}(\bm u_h-\bm{\Pi}_k^{\rm m}(\bm u,\bm u_h) ))\|_0,\label{co-02}
	\end{align}
	where $\sigma\in (1/2,1]$ is defined in Lemma \ref{embed}.
	Now by taking $\bm{\mathcal B}_k=\bm{\Pi}_k^o$ in \eqref{Bk}--\eqref{newpi2}, we use the trace inequality and an inverse inequality to get
	\begin{align}\label{newpi0-h}
	\|\bm{\Pi}_k^m(\bm u,\bm u_h)-\bm{\Pi}^o_k\bm u\|_0 &\le C \left(\|h_F^{1/2}\bm n\times[\![\bm u_h-\bm{\Pi}^o_k\bm u]\!]\|_{0,\mathcal{F}_h} + \|\bm u_h-\bm{\Pi}^o_k\bm u\|_0 \right)\le C\|\bm{\Pi}_k^o\bm u-\bm u_h\|_0.
	\end{align}
	Similarity, we can get
	\begin{align}
	\|\bm{\Pi}_{h,k}^{\rm curl,c}(\bm u_h-\bm{\Pi}_k^{\rm m}(\bm u,\bm u_h))\|_{0}&\le C\|\bm{\Pi}_k^o\bm u-\bm u_h\|_0,\label{newpi1-h}\\
	\|\nabla\times(\bm{\Pi}_{h,k}^{\rm curl,c}(\bm u_h-\bm{\Pi}_k^{\rm m}(\bm u,\bm u_h)))\|_{0}&\le C\left(\|h_F^{-1/2}\bm n\times[\![\bm u_h-\bm{\Pi}^o_k\bm u]\!]\|_{0,\mathcal{F}_h}+ \|\nabla_h\times(\bm{\Pi}^o_k\bm u-\bm u_h)\|_0\right)\nonumber\\
	&\le C\left(\|\bm\sigma-\bm\sigma_h\|_{\bm\Sigma_h}+\|\nabla_h\times(\bm{\Pi}_k^o\bm u-\bm u_h)\|_0+h^{s-1}\|\bm u\|_s\right)\nonumber\\&\le C  h^{s-1}(\|\bm{r}\|_{s-1}+\|\bm{u}\|_{s}+\|p\|_{t} ).\label{newpi2-h}
	\end{align}
	It follows from \eqref{co-02} and \eqref{newpi2-h} that
	\begin{align}
	\|\bm\Theta-( \bm{\Pi}_{h,k}^{\rm curl,c}(\bm u_h-\bm{\Pi}_k^{\rm m}(\bm u,\bm u_h) ))\|_{0}\le Ch^{s-1+\sigma}(\|\bm{r}\|_{s-1}+\|\bm{u}\|_{s}+\|p\|_{t} ).
	\label{theta-h-1}
	\end{align}	
	By \eqref{newpi1-h} and \eqref{theta-h-1}, it holds
	\begin{align}
	%\|\bm\Theta\|_{\sigma}
	%&\le C\|\bm\sigma-\bm\sigma_h\|
	%\le Ch^{s-1}(\|\bm{r}\|_{s-1}+\|\bm{u}\|_{s}+\|p\|_{t} ),\label{theta-h-2}\\
	\|\bm\Theta\|_0&\le	\|\bm\Theta-( \bm{\Pi}_{h,k}^{\rm curl,c}(\bm u_h-\bm{\Pi}_k^{\rm m}(\bm u,\bm u_h) ))\|_{0}
	+\| \bm{\Pi}_{h,k}^{\rm curl,c}(\bm u_h-\bm{\Pi}_k^{\rm m}(\bm u,\bm u_h) )\|_{0}\nonumber\\
	&\le
	C  h^{s-1+\sigma}(\|\bm{r}\|_{s-1}+\|\bm{u}\|_{s}+\|p\|_{t} )
	+C\|\bm{\Pi}_k^o\bm u-\bm u_h\|_0.\label{theta-h-2}
	\end{align}
	Similar to the proof of Lemma \ref{Jh}, we have
	\begin{align*}%\label{eq:BhJST}
	B_h(\bm{\mathcal{J}}_h\bm{\sigma}^d,\bm{\tau}_h)=-(\bm\Theta,\bm v_h)+(\Pi^{o}_kp-p_h,q_h)+E^{\mathcal{J}}_h(\bm{\sigma};\bm{\tau}_h)\qquad\forall\, \bm{\tau}_h\in\bm{\Sigma}_h.
	\end{align*}
	We take $\bm{\tau}_h=\bm{\mathcal{J}}_h\bm{\sigma}-\bm{\sigma}_h$ in last equality and use \eqref{510}, \eqref{eq:EJSJh}--\eqref{eq:EJSJhD}, \eqref{theta-h-2} to get
	\begin{align*}
	&-(\bm\Theta,\bm{\Pi}^o_k\bm u-\bm u_h)+\|\Pi^{o}_kp-p_h\|^2_0\nonumber\\
	&\qquad= B_h(\bm{\mathcal{J}}_h\bm{\sigma}^d,{\bm{\mathcal{J}}_h}\bm{\sigma}-\bm{\sigma}_h)-E^{\mathcal{J}}_h(\bm{\sigma}^d;{\bm{\mathcal{J}}_h}\bm{\sigma}-\bm{\sigma}_h) \nonumber\\
	&\qquad= E^{\mathcal{J}}_h(\bm{\sigma};\bm{\mathcal{J}}_h\bm{\sigma}^d)-E^{\mathcal{J}}_h(\bm{\sigma}^d;{\bm{\mathcal{J}}_h}\bm{\sigma}-\bm{\sigma}_h) \nonumber\\
	&\qquad\color{black}\le C h^{s-1+\beta}(\|\bm{r}\|_{s-1}+\|\bm{u}\|_{s}+\|p\|_{t} )\left(\|\bm{\Theta}\|_0+\|\Pi^{o}_kp-p_h\|_0\right)\nonumber\\
	&\qquad\le C(h^{2s-2+\sigma+\beta}+h^{2s-2+2\beta})(\|\bm{r}\|_{s-1}+\|\bm{u}\|_{s}+\|p\|_{t} )^2+\frac{1}{2}\left(\|\bm{\Pi}_k^o\bm u-\bm u_h\|^2_0+\|\Pi_k^op-p_h\|^2_0\right).
	\end{align*}
	It directly implies that
	\begin{align}\label{ppp}
	-(\bm\Theta,\bm{\Pi}^o_k\bm u-\bm u_h)+\frac{1}{2}\|\Pi^{o}_kp-p_h\|^2_0\le C(h^{2s-2+\sigma+\beta}	+h^{2s-2+2\beta})
	(\|\bm{r}\|_{s-1}+\|\bm{u}\|_{s}+\|p\|_{t} )^2+\frac{1}{2}\|\bm{\Pi}_k^o\bm u-\bm u_h\|^2_0.
	\end{align}
	We let $\sigma_h\in \mathbb P_k(\mathcal T_h)\cap H^1_0(\Omega)$ be as defined in \eqref{def2}, then we take $q_h=\widehat q_h=\sigma_h$ in \eqref{fhm3} to get
	\begin{align}\label{onn}
	-(\bm u_h,\nabla\sigma_h)=(g,\sigma_h).
	\end{align}
	We use a direct calculation to get
	\begin{align}
	&(\bm{\Pi}_k^o\bm u-\bm{\Pi}_k^{\rm m}(\bm u,\bm u_h),\bm{\Pi}_k^o\bm u-\bm u_h   )\nonumber\\
	&\qquad=
	(\bm{\Pi}_k^o\bm u-\bm{\Pi}_k^{\rm m}(\bm u,\bm u_h),\bm u-\bm u_h   )
	+(\bm{\Pi}_k^o\bm u-\bm{\Pi}_k^{\rm m}(\bm u,\bm u_h),\bm{\Pi}_k^o\bm u-\bm u   )\nonumber\\
	&\qquad=(-\nabla\sigma_h,\bm u-\bm u_h   )
	+(\bm{\Pi}_k^o\bm u-\bm{\Pi}_k^{\rm m}(\bm u,\bm u_h),\bm{\Pi}_k^o\bm u-\bm u   )
	&\text{by the definiton of $\bm\Pi^{\rm m}_k$}
	\nonumber\\
	&\qquad=(\sigma_h,\nabla\cdot\bm u   )
	+(\nabla\sigma_h,\bm u_h) 
	+(\bm{\Pi}_k^o\bm u-\bm{\Pi}_k^{\rm m}(\bm u,\bm u_h),\bm{\Pi}_k^o\bm u-\bm u   )
	&\text{by integration by parts}
	\nonumber\\
	&\qquad=(\bm{\Pi}_k^o\bm u-\bm{\Pi}_k^{\rm m}(\bm u,\bm u_h),\bm{\Pi}_k^o\bm u-\bm u   )
	&\text{by }\eqref{o1}, \eqref{onn}
	\nonumber\\
	&\qquad\le C\|\bm{\Pi}_k^o\bm u-\bm u_h\|_0\|\bm{\Pi}_k^o\bm u-\bm u \|_0
	&\text{by }\eqref{newpi0-h}. \label{4444}
	\end{align}
	We use \eqref{ic1} to get
	\begin{align}\label{iccc}
	&\|\bm{\Pi}_{h,k}^{\rm curl,c}(\bm u_h-\bm{\Pi}_k^{\rm m}(\bm u,\bm u_h) )-(\bm u_h-\bm{\Pi}_k^{\rm m}(\bm u,\bm u_h)\|_0\nonumber\\
	&\qquad\le Ch\|h_F^{-1/2}\bm n\times[\![\bm u_h-\bm{\Pi}_k^{\rm m}(\bm u,\bm u_h) ]\!]\|_{0,\mathcal F_h}\nonumber\\
	&\qquad= Ch\|h_F^{-1/2}\bm n\times[\![\bm u_h-\bm{\Pi}_k^o\bm u ]\!]\|_{0,\mathcal F_h}
	\nonumber\\
	&\qquad= 
	Ch\|h_F^{-1/2}\bm n\times[\![\bm u_h-\widehat{\bm u}_h-\bm{\Pi}_k^o\bm u +\bm{\Pi}_k^{\partial}\bm u]\!]\|_{0,\mathcal F_h}
	\nonumber\\
	&\qquad\le 
	Ch\|\bm\sigma_h-\bm{\mathcal J}_h\bm\sigma\|_{\bm{\Sigma}_h} \nonumber\\
	&\qquad\le  C  h^{s}(\|\bm{r}\|_{s-1}+\|\bm{u}\|_{s}+\|p\|_{t} ).
	\end{align}
	By using \eqref{ppp}, \eqref{theta-h-1}, \eqref{4444} and \eqref{iccc}, one can obtain
	\begin{align*}
	\|\bm{\Pi}^o_k\bm u-\bm u_h\|^2_0
	&=(\bm{\Pi}^o_k\bm u-\bm u_h,\bm{\Pi}^o_k\bm u-\bm u_h   )\nonumber\\
	&=-(\bm\Theta,\bm{\Pi}^o_k\bm u-\bm u_h   )
	+(\bm\Theta-( \bm{\Pi}_{h,k}^{\rm curl,c}(\bm u_h-\bm{\Pi}_k^{\rm m}(\bm u,\bm u_h) )),\bm{\Pi}^o_k\bm u-\bm u_h   )
	\nonumber\\
	&\quad+(\bm{\Pi}^o_k\bm u-\bm{\Pi}_k^{\rm m}(\bm u,\bm u_h),\bm{\Pi}^o_k\bm u-\bm u_h   )\nonumber\\
	&\quad+(\bm{\Pi}_{h,k}^{\rm curl,c}(\bm u_h-\bm{\Pi}_k^{\rm m}(\bm u,\bm u_h) )-(\bm u_h-\bm{\Pi}_k^{\rm m}(\bm u,\bm u_h) ),\bm{\Pi}^o_k\bm u-\bm u_h)
	\nonumber\\
	&\le
	C\left(h^{2s-2+\sigma+\beta}
	+h^{2s-2+2\beta}
	\right)
	(\|\bm{r}\|_{s-1}+\|\bm{u}\|_{s}+\|p\|_{t} )^2\nonumber\\
	&\quad+C\|\bm{\Pi}_k^o\bm u-\bm u\|^2_0
	+\frac{3}{4}
	\|\bm{\Pi}_k^o\bm u-\bm u_h\|^2_0,
	\end{align*}
	which together with \eqref{ppp} and \eqref{theta-h-2} implies \eqref{eq:PkUP}.
	\color{black}
	
\end{proof}

\subsection{Low Regularity Case}

In this subsection, we take $\bm{T}_{\Gamma,k}=\bm{\mathcal{P}}^{{\color{black}\rm div}}_{\Gamma,k}$ and $\alpha=-1$.

\subsubsection{Primary error estimates}

\begin{lemma} \label{lemma57}Let $(\bm{r},\bm{u},p)$ be the solution of \eqref{mixed},  and let $\bm{\sigma}$ and $\bm{\mathcal{I}}_h\bm{\sigma}$ be defined as
	\begin{align*}
	\bm{\sigma}&:=(\bm{r},\bm{u},\widehat{\bm{u}},p,\widehat{p}),\quad \text{ where }\, (\widehat{\bm{u}},\widehat{p}):=(\bm{u},p) \text{ on } \mathcal{F}_h,\qquad
	\bm{\mathcal{I}}_h\bm{\sigma}:=(\bm{\Pi}^{o}_m\bm{r},\bm{\mathcal{P}}^{{\color{black}\rm curl}}_k\bm{u},
	\bm{n}\times\bm{\mathcal{P}}^{{\color{black}\rm curl}}_k\bm{u}\times\bm{n}, \Pi^{o}_kp, \Pi^{\partial}_kp).
	\end{align*}
	Then we have
	\begin{equation}\label{eq:LRBhest}
	B_h(\bm{\mathcal{I}}_h\bm{\sigma},\bm{\tau}_h)=F_h(\bm{\tau}_h)+E^{\mathcal{I}}_h(\bm{\sigma};\bm{\tau}_h)\qquad \forall\,\bm{\tau}_h\in\bm{\Sigma}^{\bm{0}}_h,
	\end{equation}
	where
	\begin{align}\label{eq:DefEIhst}
	E^{\mathcal{I}}_h(\bm{\sigma};\bm{\tau}_h)&= (\nabla\times\bm{u}-\nabla\times\bm{\mathcal{P}}^{{\color{black}\rm curl}}_k\bm{u},\bm{s}_h) -\langle \bm{n}\times(\widehat{\bm{v}}_h-\bm{v}_h),\bm{\Pi}^{o}_m\bm{r}-\bm{r} \rangle_{\partial\mathcal{T}_h} \nonumber\\
	&\quad+(\bm u-\bm{\mathcal{P}}^{{\color{black}\rm curl}}_k\bm{u},\nabla_h q_h)-\langle \bm{n}\cdot(\bm{\mathcal{P}}^{{\color{black}\rm curl}}_k\bm u-\bm{u}),\widehat{q}_h -q_h \rangle_{\partial\mathcal{T}_h}
	+\langle h_F^{-1}(\Pi^{o}_kp-p),q_h-\widehat{q}_h \rangle_{\partial\mathcal{T}_h}.
	\end{align}
\end{lemma}

\begin{proof}
	By the definitions of $a_h$ and $b_h$, we get
	\begin{align*}
	&a_h(\bm{\Pi}^{o}_m\bm{r},\bm{s}_h)+b_h(\bm{\mathcal{P}}^{{\color{black}\rm curl}}_k\bm{u},\bm{n}\times\bm{\mathcal{P}}^{{\color{black}\rm curl}}_k\bm{u}\times\bm{n};\bm{s}_h)\nonumber\\
	&\qquad\quad= (\bm{\Pi}^{o}_m\bm{r},\bm{s}_h)-(\bm{\mathcal{P}}^{{\color{black}\rm curl}}_k\bm{u},\nabla_h\times\bm{s}_h)-\langle \bm{n}\times\bm{\mathcal{P}}^{{\color{black}\rm curl}}_k\bm{u},\bm{s}_h \rangle_{\partial\mathcal{T}_h} \nonumber\\
	&\qquad\quad= (\bm{r},\bm{s}_h)-(\nabla\times\bm{\mathcal{P}}^{{\color{black}\rm curl}}_k\bm{u},\bm{s}_h) \nonumber\\
	&\qquad\quad= (\nabla\times\bm{u}-\nabla\times\bm{\mathcal{P}}^{{\color{black}\rm curl}}_k\bm{u},\bm{s}_h).
	\end{align*}
	By the definitions of $b_h$, $c_h$ and $s_h^u$, we get
	\begin{align*}
	&b_h(\bm{v}_h,\widehat{\bm{v}}_h;\bm{\Pi}^{o}_m\bm{r})+c_h(\Pi^{o}_kp,\Pi^{\partial}_kp;\bm{v}_h)-s^u_h(  \bm{\mathcal{P}}^{{\color{black}\rm curl}}_k\bm{u},\bm{n}\times\bm{\mathcal{P}}^{{\color{black}\rm curl}}_k\bm{u}\times\bm{n}; \bm{v}_h,\widehat{\bm{v}}_h ) \nonumber\\
	&\qquad\quad=-(\bm{v}_h,\nabla_h\times\bm{\Pi}^{o}_m\bm{r})-\langle \bm{n}\times\widehat{\bm{v}}_h,\bm{\Pi}^{o}_m\bm{r} \rangle_{\partial\mathcal{T}_h} +(\nabla_h\cdot\bm{v}_h, \Pi^{o}_kp)-\langle \bm{n}\cdot\bm{v}_h,\Pi^{\partial}_kp \rangle_{\partial\mathcal{T}_h} \nonumber\\
	&\qquad\quad=-(\bm{v}_h,\nabla\times\bm{r})-\langle \bm{n}\times(\widehat{\bm{v}}_h-\bm{v}_h),\bm{\Pi}^{o}_m\bm{r}-\bm{r} \rangle_{\partial\mathcal{T}_h}-(\bm{v}_h, \nabla p)-\langle \bm{n}\cdot\bm{v}_h,\Pi^{\partial}_kp -p \rangle_{\partial\mathcal{T}_h} \nonumber\\
	&\qquad\quad=-(\bm{f},\bm{v}_h)-\langle \bm{n}\times(\widehat{\bm{v}}_h-\bm{v}_h),\bm{\Pi}^{o}_m\bm{r}-\bm{r} \rangle_{\partial\mathcal{T}_h}.
	\end{align*}
	By the definitions of $c_h$ and $s_h^p$, we get
	\begin{align*}
	&c_h(q_h,\widehat{q}_h;\bm{\mathcal{P}}^{{\color{black}\rm curl}}_k\bm{u})+s^p_h(\Pi^{o}_kp,\Pi^{\partial}_kp; q_h,\widehat{q}_h) \nonumber\\
	&\quad=(\nabla_h\cdot\bm{\mathcal{P}}^{{\color{black}\rm curl}}_k\bm{u}, q_h)-\langle \bm{n}\cdot\bm{\mathcal{P}}^{{\color{black}\rm curl}}_k\bm{u},\widehat{q}_h \rangle_{\partial\mathcal{T}_h}+\langle h_F^{-1}(\Pi^{o}_kp-\Pi^{\partial}_kp),q_h-\widehat{q}_h \rangle_{\partial\mathcal{T}_h} \nonumber\\
	&\quad=(\nabla\cdot\bm{u}, q_h)-(\bm{\mathcal{P}}^{{\color{black}\rm curl}}_k\bm u-\bm{u},\nabla_h q_h)-\langle \bm{n}\cdot(\bm{\mathcal{P}}^{{\color{black}\rm curl}}_k\bm u-\bm{u}),\widehat{q}_h -q_h \rangle_{\partial\mathcal{T}_h}\nonumber\\
	&\qquad\quad+\langle h_F^{-1}(\Pi^{o}_kp-\Pi^{\partial}_kp),q_h-\widehat{q}_h \rangle_{\partial\mathcal{T}_h} \nonumber\\
	&\quad=(\bm u-\bm{\mathcal{P}}^{{\color{black}\rm curl}}_k\bm{u},\nabla_h q_h)-\langle \bm{n}\cdot(\bm{\mathcal{P}}^{{\color{black}\rm curl}}_k\bm u-\bm{u}),\widehat{q}_h -q_h \rangle_{\partial\mathcal{T}_h}+\langle h_F^{-1}(\Pi^{o}_kp-\Pi^{\partial}_kp),q_h-\widehat{q}_h \rangle_{\partial\mathcal{T}_h}.
	\end{align*}
	Then the desired result \eqref{eq:LRBhest} follows immediately.
	
\end{proof}

\begin{lemma} \label{lemma58}
	Let $\bm{\sigma}$ be defined as in Lemma \ref{lemma57} and $(\bm{r},\bm{u},p)\in [H^{s}(\Omega)]^3\times[H^{s}(\Omega)]^3\times H^{s+1}(\Omega)$ with $s\in (1/2,k]$, then it holds that
	\begin{align}\label{eq:Lem510Est}
	E^{\mathcal{I}}_h(\bm{\sigma};\bm{\tau}_h)\color{black}\le C \color{black} h^{s}(\|\bm{r}\|_{s}+\|\bm{u}\|_{s}+\|p\|_{s+1} ) \|\bm{\tau}_h\|_{\bm{\Sigma}_h}.
	\end{align}
\end{lemma}
\begin{proof}
	To simplify the notation, we define
	\begin{align}
	E^{\mathcal{I}}_1&= (\nabla\times\bm{u}-\nabla\times\bm{\mathcal{P}}^{{\color{black}\rm curl}}_k\bm{u},\bm{s}_h),
	\qquad E^{\mathcal{I}}_2=-\langle \bm{n}\times(\widehat{\bm{v}}_h-\bm{v}_h),\bm{\Pi}^{o}_m\bm{r}-\bm{r} \rangle_{\partial\mathcal{T}_h}, \nonumber\\
	E^{\mathcal{I}}_3&=(\bm u-\bm{\mathcal{P}}^{{\color{black}\rm curl}}_k\bm{u},\nabla_h q_h),\qquad\qquad\quad
	E^{\mathcal{I}}_4=-\langle \bm{n}\cdot(\bm{\Pi}^{o}_k\bm u-\bm{u}),\widehat{q}_h -q_h \rangle_{\partial\mathcal{T}_h},\nonumber\\
	E^{\mathcal{I}}_5&=\langle h_F^{\alpha}(\Pi^{o}_kp-p),q_h-\widehat{q}_h \rangle_{\partial\mathcal{T}_h}.\nonumber
	\end{align}
	We then estimate $E^{\mathcal{I}}_i$ term by term as follows.
	\begin{align}
	|E^{\mathcal{I}}_1|&\color{black}\le C \color{black} h^{s}\|\bm{r}\|_{s}\|\bm{\tau}_h\|_{\bm{\Sigma}_h},\label{eq:Lem510_1}\\
	|E^{\mathcal{I}}_2|&\color{black}\le C \color{black} \sum_{T\in\mathcal{T}_h}\|\bm{n}\times(\widehat{\bm{v}}_h-\bm{v}_h)\|_{0,\partial T}h_T^{s-1/2}\|\bm{r}\|_{s,T}\le C \color{black} h^{s}\|\bm{r}\|_{s}\|\bm{\tau}_h\|_{\bm{\Sigma}_h},\\
	|E^{\mathcal{I}}_3|&\color{black}\le C \color{black} h^s\|\bm{u}\|_{s}\|\nabla_h q_h\|_0\le C \color{black} h^{s}\|\bm{u}\|_{s}\|\bm{\tau}_h\|_{\bm{\Sigma}_h},\\
	|E^{\mathcal{I}}_4|&\color{black}\le C \color{black}  \sum_{T\in\mathcal{T}_h}\|q_h-\widehat{q}_h\|_{0,\partial T}h_T^{s-1/2}\|\bm{u}\|_{s,T}\le C \color{black} h^{s}\|\bm{u}\|_{s}\|\bm{\tau}_h\|_{\bm{\Sigma}_h},\\
	|E^{\mathcal{I}}_5|&\color{black}\le C \color{black} \sum_{T\in\mathcal{T}_h}\|q_h-\widehat{q}_h\|_{0,\partial T}h_T^{s-1/2}\|p\|_{s+1,T}\le C \color{black} h^{s}\|p\|_{s+1}\|\bm{\tau}_h\|_{\bm{\Sigma}_h}.\label{eq:Lem510_5}
	\end{align}
	In view of \eqref{eq:DefEIhst}, \eqref{eq:Lem510_1}--\eqref{eq:Lem510_5}, we get the desired estimation \eqref{eq:Lem510Est}.
\end{proof}

\begin{theorem} \label{l53344} Let $\bm{\sigma}$ be defined as in Lemma \ref{lemma57} and $(\bm{r},\bm{u},p)\in [H^{s}(\Omega)]^3\times[H^{s}(\Omega)]^3\times H^{s+1}(\Omega)$ with $s\in (1/2,k]$, 
	let $\bm\sigma_h$ be the solution of \eqref{Bh_HDG},	
	then there holds
	\begin{eqnarray}\label{es-main}
	\|\bm{\sigma}_h-\bm{\sigma}\|_{\bm{\Sigma}_h}\color{black}\le C \color{black} h^{s}(\|\bm{r}\|_{s}+\|\bm{u}\|_{s}+\|p\|_{s+1} ).
	\end{eqnarray}
\end{theorem}
\begin{proof} By Lemma \ref{lemma57} and Lemma \ref{lemma58} we get
	\begin{align*}
	\|\bm{\sigma}_h-\bm{\mathcal{I}}_h\bm{\sigma}\|_{\bm{\Sigma}_h}
	\color{black}\le C \color{black}\sup_{\bm 0\neq\bm{\tau}_h\in \bm{\Sigma}_h^{\bm 0}}\frac{B_h( \bm{\sigma}_h-\bm{\mathcal{I}}_h\bm{\sigma},\bm{\tau}_h  )}{\|\bm{\tau}_h\|_{\bm{\Sigma}_h}} =\sup_{\bm 0\neq\bm{\tau}_h\in \bm{\Sigma}_h^{\bm 0}}\frac{E^{\mathcal{I}}_h( \bm{\sigma}_h,\bm{\tau}_h  )}{\|\bm{\tau}_h\|_{\bm{\Sigma}_h}} \le C \color{black} h^{s}(\|\bm{r}\|_{s}+\|\bm{u}\|_{s}+\|p\|_{s+1} ),
	\end{align*}
	which combines with a triangle inequality to get our desired result.
\end{proof}

\subsubsection{$L^2$ error estimates}
Assume $\bm\Theta\in \bm{H}({\color{black}\rm div} ^0;\Omega)$.
We also introduce the following dual problem:
\begin{eqnarray}
\left\{
\begin{aligned}
\bm{r}^d - \nabla\times\bm{u}^d&=0,&\text{in }\Omega,\\
\nabla\times \bm{r}^d
+\nabla p^d&=\bm\Theta, &\text{in }\Omega, \label{dual}\\
\nabla\cdot\bm{u}^d&=\Lambda,  &\text{in }\Omega,  \\
\bm{n}\times \bm{u}^d &=\bm{0}, &\text{on }\Gamma,\\
p^d&=0, &\text{on }\Gamma.
\end{aligned}
\right.
\end{eqnarray}
Using Theorem \ref{reg}, there exists a regularity index $s_0\in\left(1/2,1\right]$ such that
\begin{align}\label{reg2}
\|\bm{r}^d\|_{s_0}+\|\bm{u}^d\|_{s_0}\color{black}\le C\left( \color{black} \|\bm\Theta\|_0+\|\Lambda\|_0\right),\text{ and } p^d=0.
\end{align}

\begin{lemma} Let  $(\bm r,\bm u, p)$ and $(\bm r^d,\bm u^d,p^d)$ be the solution of \eqref{mixed} and \eqref{dual}, respectively, let $\bm\sigma_h$ be the solution of \eqref{Bh_HDG}, then we have
	\begin{align}\label{tilde}
	B_h(\bm{\sigma}^d-{\bm{\mathcal{I}}}_h\bm{\sigma}^d,\bm{\sigma}-\bm{\sigma}_h)\le   Ch^{s_0}(\|\bm r^d\|_{s_0}+\|\bm u^d\|_{s_0})\|\bm \sigma-\bm\sigma_h\|_{\bm\Sigma_h}.
	\end{align}
	Here
	$\bm{\sigma}^d:=(\bm{r}^d,\bm{u}^d,\widehat{\bm{u}}^d,p^d,\widehat{p}^d)$, with $(\widehat{\bm{u}}^d,\widehat{p}^d)=(\bm{u}^d,p^d)$ on $\mathcal{F}_h$.
\end{lemma}
\begin{proof}
	\color{black}
	From the definition of  $B_h$ and integration by parts we have
	\begin{align*}
	&B_h(\bm{\sigma}^d-{\bm{\mathcal{I}}}_h\bm{\sigma}^d,\bm{\sigma}-\bm{\sigma}_h)\nonumber\\
	&\qquad=(\bm r^d-\bm{\Pi}_m^{o}\bm r^d,\bm r-\bm r_h)
	-(\nabla\times(\bm u^d-\bm{\mathcal P}_k^{\rm curl}\bm u^d),\bm r-\bm r_h)\nonumber\\
	&\qquad\quad-
	(\bm r^d-\bm{\Pi}_m^{o}\bm r^d,\nabla_h\times(\bm u-\bm u_h))
	-\langle\bm r^d-\bm{\Pi}_m^{o}\bm r^d,\bm n\times(\widehat{\bm u}-\widehat{\bm u}_h-(\bm u-\bm u_h)) \rangle_{\partial\mathcal T_h} \nonumber\\
	&\qquad\quad-(\bm u^d-\bm{\mathcal P}_k^{\rm curl}\bm u^d,\nabla_h(p-p_h))
	+\langle \bm n\cdot(\bm u^d-\bm{\mathcal P}_k^{\rm curl}\bm u^d),(p-p_h)-(\widehat p-\widehat p_h) \rangle_{\partial\mathcal T_h} \nonumber\\
	&\qquad\le Ch^{s_0}(\|\bm r^d\|_{s_0}+\|\bm u^d\|_{s_0})\|\bm \sigma-\bm\sigma_h\|_{\bm\Sigma_h},
	\end{align*}
	where we have used the fact $p^d=0$.
\end{proof}

\begin{theorem}\label{L^222}
	Let $(\bm{r},\bm{u},p)$ be the solution of \eqref{mix0source} and $(\bm{r},\bm{u},p)\in [H^{s}(\Omega)]^3\times[H^{s}(\Omega)]^3\times H^{s+1}(\Omega)$ with $s\in (1/2,k]$, 
	let $(\bm r_h,\bm u_h,\widehat{\bm u}_h,p_h,\widehat p_h)$ be the solution of \eqref{fem},
	if $\bm g_T\in [ H^{s_g}(\Gamma)]^3$ with $s_g>1/2$, then there holds
	\begin{align}\label{eq:Thm513Est}
	\|\bm{u}-\bm{u}_h\|_0+\|p-p_h\|_0&\le C\left(
	h^{s+(s_0+\sigma)/2}+h^{s+s_0}
	\right)
	(\|\bm{r}\|_s+\|\bm{u}\|_{s}+\|p\|_{s+1})\nonumber\\
	&\quad+Ch^{\min(s_0,k-1)+s_g-1/2}\|\bm g_T\|_{s_g,\Gamma}+C\|\bm u-\bm{\mathcal P}_k^{\rm curl}\bm u\|_0,
	\end{align}
	where $\sigma\in (1/2,1]$ is defined in Lemma \ref{embed}, $s_0\in (1/2,1]$ is defined in \eqref{reg2}; otherwise, if  $\bm g_T\not\in [ H^{s_g}(\Gamma)]^3$ with any $s_g>1/2$, it holds
	\begin{align}\label{es2}
	\|\bm{u}-\bm{u}_h\|_0+\|p-p_h\|_0&\le C  \left(
	h^{s+(\min(s_0,k-1)+\sigma)/2}
	+h^{s+\min(s_0,k-1)}
	\right)
	(\|\bm{r}\|_s+\|\bm{u}\|_s+\|p\|_{s+1})\nonumber\\
	&\quad+C\|\bm u-\bm{\mathcal P}_k^{\rm curl}\bm u\|_0.
	\end{align}
\end{theorem}
\begin{proof}	\color{black}
	We first prove \eqref{eq:Thm513Est}. We take
	$\Lambda=p-p_h$ in \eqref{dual} and let $\bm\Theta\in \bm H({\rm curl};\Omega)\cap \bm H({\rm div};\Omega)$ be the solution of
	\begin{subequations}
		\begin{align*}
		\nabla\times\bm\Theta&=\nabla\times ( \bm{\Pi}_{h,k}^{\rm curl,c}(\bm u_h-\bm{\Pi}_k^{\rm m}(\bm u,\bm u_h)))&\text{in }\Omega,\\
		\nabla\cdot\bm\Theta&=0&\text{in }\Omega,\\
		\bm n\times\bm\Theta&=\bm 0&\text{on }\Gamma.
		\end{align*}
	\end{subequations}
	Due to \eqref{oror} and the result in \cite[Lemma 4.5]{MR2009375} one has
	\begin{align}
	%\|\bm\Theta\|_{\sigma}&\le C\|\nabla\times ( \bm{\Pi}_{h,k}^{\rm curl,c}(\bm u_h-\bm{\Pi}_k^{\rm m}(\bm u,\bm u_h) ))\|_0,\label{co-1}\\
	\|\bm\Theta-( \bm{\Pi}_{h,k}^{\rm curl,c}(\bm u_h-\bm{\Pi}_k^{\rm m}(\bm u,\bm u_h) ))\|_{0}&\le Ch^{\sigma}\|\nabla\times ( \bm{\Pi}_{h,k}^{\rm curl,c}(\bm u_h-\bm{\Pi}_k^{\rm m}(\bm u,\bm u_h) ))\|_0,\label{co-2}
	\end{align}
	where $\sigma\in (1/2,1]$ is defined in Lemma \ref{embed}. Now taking $\bm{\mathcal B}_k=\bm{\mathcal P}^{\rm curl}_k$ in \eqref{Bk}--\eqref{newpi2}, we obtain the following estimates,
	\begin{align}\label{newpi0-l}
	\|\bm{\mathcal P}_k^{\rm curl}\bm u-\bm{\Pi}_k^{\rm m}(\bm u,\bm u_h)\|_0
	&\le\left(\|h_F^{1/2}\bm n\times[\![\bm u_h-\bm{\mathcal P}^{\rm curl}_k\bm u]\!]\|_{0,\mathcal{F}_h}+ \|\bm u_h-\bm{\mathcal P}^{\rm curl}_k\bm u\|_0\right) \le C\|\bm{\mathcal P}^{\rm curl}_k\bm u-\bm u_h\|_0.
	\end{align}
	Similarity, we can get
	\begin{align}
	\|\bm{\Pi}_{h,k}^{\rm curl,c}(\bm u_h-\bm{\Pi}_k^{\rm m}(\bm u,\bm u_h))\|_{0}
	&\le C\|\bm{\mathcal P}^{\rm curl}_k\bm u-\bm u_h\|_0,\label{newpi1-l}\\
	\|\nabla\times(\bm{\Pi}_{h,k}^{\rm curl,c}(\bm u_h-\bm{\Pi}_k^{\rm m}(\bm u,\bm u_h)))\|_{0}
	&\le C \left( \|h_F^{-1/2}\bm n\times[\![\bm u_h-\bm{\mathcal P}^{\rm curl}_k\bm u]\!]\|_{0,\mathcal{F}_h} + \|\nabla_h\times(\bm{\mathcal P}^{\rm curl}_k\bm u-\bm u_h)\|_0 \right)\nonumber\\
	&\le C\left(\|\bm\sigma-\bm\sigma_h\|_{\bm\Sigma_h}+\|\nabla_h\times(\bm{\mathcal P}^{\rm curl}_k\bm u-\bm u_h)\|_0\right)\nonumber\\
	&\le C\left(\|\bm\sigma-\bm\sigma_h\|_{\bm\Sigma_h}+h^{s}\|\bm r\|_s\right)\nonumber\\
	&\le Ch^{s}(\|\bm{r}\|_{s}+\|\bm{u}\|_{s}+\|p\|_{s+1} ).\label{newpi2-l}
	\end{align}
	It then follows from \eqref{co-2}  and \eqref{newpi2-l} that
	\begin{align}
	\|\bm\Theta-( \bm{\Pi}_{h,k}^{\rm curl,c}(\bm u_h-\bm{\Pi}_k^{\rm m}(\bm u,\bm u_h) ))\|_{0}&\le Ch^{s+\sigma}(\|\bm{r}\|_{s}+\|\bm{u}\|_{s}+\|p\|_{s+1} ).\label{theta-l-1}
	\end{align}	
	Follows from the above estimates inequality, it holds that
	\begin{align}
	\|\bm\Theta\|_0&\le
	\|\bm\Theta-( \bm{\Pi}_{h,k}^{\rm curl,c}(\bm u_h-\bm{\Pi}_k^{\rm m}(\bm u,\bm u_h) ))\|_{0}
	+\| \bm{\Pi}_{h,k}^{\rm curl,c}(\bm u_h-\bm{\Pi}_k^{\rm m}(\bm u,\bm u_h) )\|_{0}\nonumber\\
	&\le C h^{s+\sigma}(\|\bm{r}\|_{s}+\|\bm{u}\|_{s}+\|p\|_{s+1} )+C\|\bm{\mathcal P}_k^{\rm curl}\bm u-\bm u_h\|_0.\label{theta-l-2}
	\end{align}
	From direct calculations, we have
	\begin{align}%\label{eq:BhJST}
	B_h(\bm{\sigma},\bm{\tau})&=-(\bm f,\bm v)
	+(g,q)
	\qquad\forall\, \bm{\tau}\in\bm{\Sigma}^{\bm 0},\label{b1}\\
	B_h(\bm{\sigma}^d,\bm{\tau})&=-(\bm\Theta,\bm v)
	+(p-p_h,q)
	-\langle\bm r^d,\bm n\times\widehat{\bm v}_h \rangle_{\Gamma}
	\qquad\forall\, \bm{\tau}\in\bm{\Sigma},\label{b2}
	\end{align}
	where $\bm\tau=(\bm s,\bm v,\widehat{\bm v}, q,\widehat q)$, $\bm\Sigma^{\bm 0}= \bm S\times\bm V\times\widehat{\bm V}^{\bm 0}\times Q\times Q^0$, 
	$\bm\Sigma= \bm S\times\bm V\times\widehat{\bm V}\times Q\times Q^0$, $\bm S$ and $\bm V$ are broken $\bm H({\rm curl})$ spaces, $Q$ is broken $H^1$ space, $\widehat{\bm V}$, $\widehat{\bm V}^{\bm 0}$ and $Q^0$ are defined as
	\begin{align*}
	\widehat{\bm V}&:=[L^2(\mathcal F_h)]^3,\\
	\widehat{\bm V}^{\bm 0}&:=\{\widehat{\bm v}\in [L^2(\mathcal F_h)]^3:\bm n\times\widehat{\bm v}|_{\Gamma}=\bm 0\},\\
	\widehat{Q}^{0}&:=\{\widehat q\in L^2(\mathcal F_h):\widehat q|_{\Gamma}= 0\}.
	\end{align*}
	By the fact $\bm n\times\widehat{\bm u}|_{\Gamma}=\bm n\times\bm u|_{\Gamma}=\bm g_T$ and $\bm n\times\widehat{\bm u}_h|_{\Gamma}=\bm{\mathcal P}_{\Gamma,k}^{\rm div}\bm g_T$, we can get
	\begin{align}\label{g_T}
	\langle\bm r^d,\bm n\times(\widehat{\bm u}-\widehat{\bm u}_h) \rangle_{\Gamma}&=\langle\bm r^d,\bm g_T-\bm{\mathcal P}_{\Gamma,k}^{\rm div}\bm g_T  \rangle_{\Gamma}\nonumber\\
	& =\langle\bm r^d-\min(1,k-1)\bm{\Pi}_0^o\bm r^d,\bm g_T-\bm{\mathcal P}_{\Gamma,k}^{\rm div}\bm g_T  \rangle_{\Gamma}\nonumber\\
	&\le Ch^{\min(s_0,k-1)+s_g-1/2}\|\bm r^d\|_{s_0}\|\bm g_T\|_{s_g,\Gamma}.
	\end{align}
	We take $\bm{\tau}=\bm{\sigma}-\bm{\sigma}_h\in \bm\Sigma$ in \eqref{b2}, to get
	\begin{align*}
	&-(\bm\Theta,\bm u-\bm u_h)+\|p-p_h\|^2_0\nonumber\\
	&\qquad= B_h(\bm{\sigma}^d,\bm{\sigma}-\bm{\sigma}_h)+\langle\bm r^d,\bm n\times(\widehat{\bm u}-\widehat{\bm u}_h)\rangle_{\Gamma}\nonumber\\
	&\qquad= 
	B_h(\bm{\sigma}^d-{\bm{\mathcal{I}}}_h\bm{\sigma}^d,\bm{\sigma}-\bm{\sigma}_h)
	+\langle\bm r^d,\bm n\times(\widehat{\bm u}-\widehat{\bm u}_h)\rangle_{\Gamma}
	&\text{by }\eqref{b1},\eqref{Bh_HDG}\nonumber\\
	&\qquad\le  Ch^{s_0}(\|\bm r^d\|_{s_0}+\|\bm u^d\|_{s_0})\|\bm \sigma-\bm\sigma_h\|_{\bm\Sigma}+Ch^{\min(s_0,k-1)+s_g-1/2}\|\bm g_T\|_{s_g,\Gamma} &\text{by }\eqref{tilde},\eqref{g_T} \nonumber\\
	&\qquad\color{black}\le C h^{s+s_0}\left(\|\bm{\Theta}\|_0+\|p-p_h\|_0\right)(\|\bm{r}\|_{s}+\|\bm u\|_s)+Ch^{\min(s_0,k-1)+s_g-1/2}\|\bm g_T\|_{s_g,\Gamma}&\text{by }\eqref{es-main},\eqref{reg2}\nonumber\\
	&\qquad\le C\left(
	h^{2s+s_0+\sigma}+h^{2s+2s_0}
	\right)(\|\bm{r}\|_{s}+\|\bm{u}\|_{s}+\|p\|_{s+1} )^2\nonumber\\
	&\qquad\quad+\frac{1}{2}\left(\|\bm u-\bm u_h\|^2_0+\|p-p_h\|^2_0\right)
	+Ch^{\min(s_0,k-1)+s_g-1/2}\|\bm g_T\|_{s_g,\Gamma}+C\|\bm u-\bm{\mathcal P}_k^{\rm curl}\bm u\|_0 &\text{by }\eqref{theta-l-2}.
	\end{align*}
	It directly implies that
	\begin{align}\label{pppp}
	&-(\bm\Theta,\bm u-\bm u_h)+\frac{1}{2}\|p-p_h\|^2_0\nonumber\\
	&\qquad\le C\left(
	h^{2s+s_0+\sigma}+h^{2s+2s_0}
	\right)(\|\bm{r}\|_{s}+\|\bm{u}\|_{s} +\|p\|_{s+1})^2+\frac{1}{2}\|\bm u-\bm u_h\|^2_0\nonumber\\
	&\qquad\quad
	+Ch^{\min(s_0,k-1)+s_g-1/2}\|\bm g_T\|_{s_g,\Gamma}
	+C\|\bm u-\bm{\mathcal P}_k^{\rm curl}\bm u\|_0 .
	\end{align}
	We let $\sigma_h\in  \mathbb P_k(\mathcal T_h)\cap H^1_0(\Omega)$ be as defined in \eqref{def2}, then we take $q_h=\widehat q_h=\sigma_h$ in \eqref{fhm3} to get
	\begin{align}\label{on}
	-(\bm u_h,\nabla\sigma_h)=(g,\sigma_h).
	\end{align}
	We use a direct calculation to get
	\begin{align}
	&(\bm u-\bm{\Pi}_k^{\rm m}(\bm u,\bm u_h),\bm u-\bm u_h   )\nonumber\\
	&\qquad=
	(\bm u-\bm{\mathcal P}_k^{\rm curl}\bm u,\bm u-\bm u_h)
	+
	(-\nabla\sigma_h,\bm u-\bm u_h   )
	&\text{by the definiton of $\bm\Pi^{\rm m}_k$}
	\nonumber\\
	&\qquad=
	(\bm u-\bm{\mathcal P}_k^{\rm curl}\bm u,\bm u-\bm u_h)
	+
	(\sigma_h,\nabla\cdot\bm u   )
	+(\nabla\sigma_h,\bm u_h) 
	&\text{by integration by parts}
	\nonumber\\
	&\qquad=(\bm u-\bm{\mathcal P}_k^{\rm curl}\bm u,\bm u-\bm u_h)
	&\text{by }\eqref{o1}, \eqref{on} \label{444}.
	\end{align}
	We use \eqref{ic1} to get
	\begin{align}\label{icc}
	&\|\bm{\Pi}_{h,k}^{\rm curl,c}(\bm u_h-\bm{\Pi}_k^{\rm m}(\bm u,\bm u_h) )-(\bm u_h-\bm{\Pi}_k^{\rm m}(\bm u,\bm u_h)\|_0\nonumber\\
	&\qquad\le Ch\|h_F^{-1/2}\bm n\times[\![\bm u_h-\bm{\Pi}_k^{\rm m}(\bm u,\bm u_h) ]\!]\|_{0,\mathcal F_h}\nonumber\\
	&\qquad= Ch\|h_F^{-1/2}\bm n\times[\![\bm u_h ]\!]\|_{0,\mathcal F_h^I}
	+
	Ch\|h_F^{-1/2}\bm n\times[\![\bm u_h-\bm{\mathcal P}_k^{\rm curl}\bm u ]\!]\|_{0,\mathcal F_h^B}
	\nonumber\\
	&\qquad= Ch\|h_F^{-1/2}\bm n\times[\![\bm u_h-\widehat{\bm u}_h ]\!]\|_{0,\mathcal F_h^I}
	+
	Ch\|h_F^{-1/2}\bm n\times[\![\bm u_h-\widehat{\bm u}_h]\!]\|_{0,\mathcal F_h^B}
	\nonumber\\
	&\qquad\le 
	Ch\|\bm\sigma_h-\bm{\mathcal I}_h\bm\sigma\|_{\bm{\Sigma}_h} \nonumber\\
	&\qquad\le  C  h^{s+1}(\|\bm{r}\|_{s}+\|\bm{u}\|_{s}+\|p\|_{s+1} ).
	\end{align}
	By using \eqref{pppp}, \eqref{theta-l-1}, \eqref{444}, and \eqref{icc}, one can obtain
	\begin{align*}
	\|\bm u-\bm u_h\|^2_0
	&=(\bm u-\bm u_h,\bm u-\bm u_h   )\nonumber\\
	&=-(\bm\Theta,\bm u-\bm u_h   )+(\bm\Theta-( \bm{\Pi}_{h,k}^{\rm curl,c}(\bm u_h-\bm{\Pi}_k^{\rm m}(\bm u,\bm u_h) )),\bm u-\bm u_h   )\nonumber\\
	&\quad+(\bm u-\bm{\Pi}_k^{\rm m}(\bm u,\bm u_h),\bm u-\bm u_h   )\nonumber\\
	&\quad+(\bm{\Pi}_{h,k}^{\rm curl,c}(\bm u_h-\bm{\Pi}_k^{\rm m}(\bm u,\bm u_h) )-(\bm u_h-\bm{\Pi}_k^{\rm m}(\bm u,\bm u_h) )
	,\bm u-\bm u_h)\nonumber\\
	&\le C\left(
	h^{2s+s_0+\sigma}+h^{2s+2s_0}
	\right)(\|\bm{r}\|_{s}+\|\bm{u}\|_{s}+\|p\|_{s+1} )^2\nonumber\\
	&\quad+
	Ch^{\min(s_0,k-1)+s_g-1/2}\|\bm g_T\|_{s_g,\Gamma}
	+\frac{1}{2}\|\bm u-\bm u_h\|^2_0+C\|\bm u-\bm{\mathcal P}_k^{\rm curl}\bm u\|_0,
	\end{align*}
	which together with \eqref{pppp} and \eqref{theta-l-2} implies \eqref{eq:Thm513Est}. 
	The proof for \eqref{es2} is followed by the similar steps in \Cref{l53344} and by taking $\bm{\mathcal B}_k=\bm{\mathcal P}^{\rm curl}_k$ in \eqref{Bk}--\eqref{newpi2}, thus we omit it.
	\color{black}
	
\end{proof}

\section{Numerical experiments}\label{sec:NE}
In this section, we simulate numerical experiments to test the performance of HDG scheme \eqref{fem}. All numerical tests are programmed in C++. The interior unknowns $\bm{r}_h$, ${\bm{u}_h}$ and ${p_h}$ are eliminated locally, and hence the resulting system consists of only the unknowns $\widehat{\bm{u}}_h$ and $\widehat{p}_h$. The quasi-uniform simplex meshes are used for all numerical examples. GMRES and SparseLU are used as the linear system solvers.
\subsection{Smooth solution}
In the first experiment, we take $\Omega=[0,1]^3$. The functions $\bm{f}$, $g$ and $\bm{g}_T$ in \eqref{source} are chosen according to the following exact solutions with high regularities.
\begin{align*}
u_1=\sin(\pi y)\sin(\pi z),\qquad
u_2=\sin(\pi y)\sin(\pi z),\qquad
u_3=\sin(\pi y)\sin(\pi z),\\
p=\sin(\pi x)\sin(\pi y)\sin(\pi z)/\pi^2.\nonumber
\end{align*}

\begin{table}[!h]
	\small \label{com}
	\caption{\label{tab1}Results for smooth solutions}
	\centering
	\vspace{0.1in}
	\begin{tabular}{c|c|c|c|c|c|c|c|c}
		\Xhline{1pt}
		\multirow{2}{*}{$k$, $m$ and $\alpha$ } &
		
		\multirow{2}{*}{$h^{-1}$} &
		
		\multicolumn{2}{c|}{${\|\bm{r}-\bm{r}_h\|_0}$} &
		\multicolumn{2}{c|}{${\|\bm{u}-\bm{u}_h\|_0}$} &
		
		\multicolumn{2}{c| }{${\|p-p_{h}\|_0}$} &
		\multirow{2}{*}{DOF}\\
		
		\cline{3-8}
		
		&  &Error &Rate  &Error &Rate  &Error &Rate  \\
		\hline

		$k=1$&2	&1.75E+00	&	    &5.41E-01	&	    &1.07E-01	&	    	&1080\\
		$m=0$&4	&9.32E-01	&0.91 	&2.06E-01	&1.39 	&2.86E-02	&1.90 		&7776\\
		$\alpha=-1$&8	&4.74E-01	&0.98 	&6.50E-02	&1.66 	&5.07E-03	&2.50 		&58752\\
		&16	&2.38E-01	&0.99 	&1.78E-02	&1.87 	&7.24E-04	&2.81 		&456192\\
		\cline{1-9}
		$k=1$&2	&1.74E+00	&	    &4.07E-01	&	    &2.35E-01	&	    	&1080\\
		$m=0$&4	&9.31E-01	&0.91 	&1.12E-01	&1.87 	&1.46E-01	&0.69 		&7776\\
		$\alpha=1$&8	&4.73E-01	&0.98 	&2.80E-02	&2.00 	&7.76E-02	&0.91 		&58752\\
		&16	&2.38E-01	&0.99 	&6.90E-03	&2.02 	&3.96E-02	&0.97 		&456192\\
		
		\cline{1-9}
		
		$k=2$&2	&1.12E-01	&	    &3.21E-02	&	    &1.35E-02	&	    	&2160\\
		$m=2$&4	&1.72E-02	&2.70 	&4.51E-03	&2.83 	&1.06E-03	&3.67 		&15552\\
		$\alpha=-1$&8	&3.04E-03	&2.50 	&6.23E-04	&2.86 	&7.50E-05	&3.82 		&117504\\
		
		\cline{1-9}
		$k=2$&2	&1.21E-01	&	    &2.67E-02	&	    &2.89E-02	&	    	&2160\\
		$m=2$&4	&2.16E-02	&2.49 	&3.13E-03	&3.09 	&5.84E-03	&2.31 		&15552\\
		$\alpha=1$&8	&4.50E-03	&2.26 	&3.81E-04	&3.04 	&1.32E-03	&2.15 		&117504\\
		
		\Xhline{1pt}
	\end{tabular}
\end{table}

From the numerical results in Table~\ref{tab1}, we observe that optimal convergence rates can be obtained for $\bm{r}$, $\bm{u}$ and $p$, confirming the theoretical results.
Moreover, it can be observed that the convergence rate for $p$ as $\alpha=-1$ is better than that as $\alpha=1$; the convergence rate for $u$ as $\alpha=1$ is better than that as $\alpha=-1$.

\subsection{Singular solution on $L$-shaped domain}

In the second experiment, we consider an $L$-shaped domain $\Omega=[-1,1]^3/(-1,0)\times(-1,0)\times(-1,1)$. The functions $\bm{f}$, $g$ and $\bm{g}_T$ in \eqref{source} are chosen according to the following exact (singular) solutions
\begin{align*}
u_1=tr^{t-1}\sin[(t-1)\theta],\quad
u_2=tr^{t-1}\cos[(t-1)\theta],\quad
u_3=0,\qquad
\bm{r}=\bm{0},\qquad
p=0.\nonumber
\end{align*}

By taking $t=\frac{2}{3}$ and $\frac{4}{3}$, we have $\bm{u}\in [H^{\frac{2}{3}-\epsilon}(\Omega)]^3$ and
$\bm{u}\in [H^{\frac{4}{3}-\epsilon}(\Omega)]^3$, respectively, for arbitrary $\epsilon>0$. The results for $k=m=1$, $\alpha=-1$ are presented in Table~\ref{L}.

\begin{table}[!h]
	\small\label{L}
	\caption{\label{tab2}Results on $L$-shaped domain with singular solutions}
	\centering
	\vspace{0.1in}
	\begin{tabular}{c|c|c|c|c|c|c|c|c}
		\Xhline{1pt}
		\multirow{2}{*}{$t$} &
		\multirow{2}{*}{$h^{-1}$} &
		
		\multicolumn{2}{c|}{$\|\bm{r}-\bm{r}_h\|_0$} &
		\multicolumn{2}{c|}{$\|\bm{u}-\bm{u}_h\|_0$} &
		\multicolumn{2}{c| }{$\|p-p_{h}\|_0$} &
		\multirow{2}{*}{DOF}\\
		\cline{3-8}
		
		&  &Error &Rate  &Error &Rate  &Error &Rate  \\
		\hline
		&2	&1.56E-01	&	    &2.40E-01	&	    &1.35E-01	&	    	&846\\
		&4	&9.73E-02	&0.69 	&1.77E-01	&0.44 	&8.00E-02	&0.76 		&5976\\
		$\frac{2}{3}$&8	&4.72E-02	&1.04 	&1.26E-01	&0.49 	&3.64E-02	&1.14 		&44640\\
		&16	&1.98E-02	&1.26 	&8.31E-02	&0.60 	&1.42E-02	&1.36 		&344448\\
		&32	&1.16E-02	&0.77 	&5.36E-02	&0.63 	&5.63E-03	&1.33 		&2704896\\
		\cline{1-9}
		
		&2	&7.19E-02	&	    &1.06E-01	&	    &7.49E-02	&	    	&846\\
		&4	&2.88E-02	&1.32 	&5.17E-02	&1.04 	&2.83E-02	&1.40 		&5976\\
		$\frac{4}{3}$&8	&8.73E-03	&1.72 	&2.40E-02	&1.11 	&7.96E-03	&1.83 		&44640\\
		&16	&2.70E-03	&1.70 	&1.06E-02	&1.17 	&1.82E-03	&2.13 		&344448\\
		&32	&5.01E-04	&2.43 	&4.12E-03	&1.37 	&4.23E-04	&2.10 		&2704896\\
		
		\Xhline{1pt}
	\end{tabular}
\end{table}

It can be observed from Table \ref{tab2} that the convergence rates of $\|\bm{u}-\bm{u}_h\|_0$ are approaching to $\frac{2}{3}$ and $\frac{4}{3}$ in the cases when $t$ equals to $\frac{2}{3}$ and $\frac{4}{3}$, respectively. It confirms the theoretical result in Theorem~\ref{L^222}. Moreover, when $t=\frac{4}{3}$, the convergence rates of $\|\bm{r}-\bm{r}_h\|_0$ and $\|p-p_h\|_0$ goes to $2$, which are better than we expected.
When $t=\frac{2}{3}$,  the convergence rates of $\|p-p_h\|_0$ tends to be twice of $\frac{2}{3}$, while the convergence rates of $\|\bm{r}-\bm{r}_h\|_0$ is around $\frac{2}{3}$.

Next we compare the errors obtained by using $\bm{\mathcal{P}}^{{\color{black}\rm div}}_{\Gamma,k}\bm{g}_T$ and $\bm{\Pi}^{\partial}_k\bm{g}_T$ to approximate $\bm{g}_T$ on $\Gamma$. One can see from Table~\ref{camp3} and Table~\ref{camp4} that, for $t=0.55$, $k=1$, $m=0 \,\,\text{or} \,\,1$, the convergence rate of $\|\bm{r}-\bm{r}_h\|_0$ deteriorates when $\bm{\Pi}^{\partial}_k\bm{g}_T$ is used. It also illustrates that the projection $\bm{\mathcal{P}}^{{\color{black}\rm div}}_{\Gamma,k}$ has an effect on enhancing accuracy in the computation.

\begin{table}[!h]
	\small
	\caption{\label{camp3}Comparison of $\bm{\mathcal{P}}^{{\color{black}\rm div}}_{\Gamma,k}\bm{g}_T$ and $\bm{\Pi}^{\partial}_k\bm{g}_T$ for $t=0.55$, $k=1$ and $m=0$}
	\centering
	
	\begin{tabular}{c|c|c|c|c|c|c|c|c}
		\Xhline{1pt}
		\multirow{1}{*}{ } &
		
		\multirow{2}{*}{$h^{-1}$} &
		
		\multicolumn{2}{c|}{$\|\bm{r}-\bm{r}_h\|_0$} &
		\multicolumn{2}{c|}{$\|\bm{u}-\bm{u}_h\|_0$} &
		\multicolumn{2}{c| }{$\|p-p_{h}\|_0$} &
		\multirow{2}{*}{DOF} \\
		\cline{3-8}
		
		&  &Error &Rate  &Error &Rate  &Error &Rate  \\
		\hline
		
		&2	&1.44E-01	&	    &3.35E-01	&	    &1.43E-01	&	    	&846\\
		&4	&1.09E-01	&0.40 	&2.59E-01	&0.37 	&1.01E-01	&0.50 		&5976\\
		$\bm{\mathcal{P}}^{{\color{black}\rm div}}_{\Gamma,k}\bm{g}_T$&8	&6.38E-02	&0.78 	&2.04E-01	&0.35 	&5.33E-02	&0.92 		&44640\\
		&16	&3.30E-02	&0.95 	&1.49E-01	&0.45 	&2.32E-02	&1.20 		&344448\\
		&32	&2.01E-02	&0.71 	&1.04E-01	&0.52 	&9.40E-03	&1.30 		&2704896\\
		
		\cline{1-9}
		
		&2	&1.64E-01	&	    &3.05E-01	&	    &1.30E-01	&	    	&846\\
		&4	&1.32E-01	&0.31 	&2.39E-01	&0.35 	&9.85E-02	&0.40 		&5976\\
		$\bm{\Pi}^{\partial}_k\bm{g}_T$&8	&9.47E-02	&0.48 	&1.92E-01	&0.32 	&5.22E-02	&0.92 		&44640\\
		&16	&8.26E-02	&0.20 	&1.42E-01	&0.43 	&2.26E-02	&1.20 		&344448\\
		&32	&9.29E-02	&-0.17 	&9.97E-02	&0.51 	&9.16E-03	&1.31 		&2704896\\
		
		\Xhline{1pt}
	\end{tabular}
\end{table}

\begin{table}[!h]
	\small
	\caption{\label{camp4}Comparison of $\bm{\mathcal{P}}^{{\color{black}\rm div}}_{\Gamma,k}\bm{g}_T$ and $\bm{\Pi}^{\partial}_k\bm{g}_T$ for $t=0.55$, $k=1$ and $m=1$}
	\centering
	
	\begin{tabular}{c|c|c|c|c|c|c|c|c}
		\Xhline{1pt}
		\multirow{1}{*}{ } &
		
		\multirow{2}{*}{$h^{-1}$} &
		
		\multicolumn{2}{c|}{$\|\bm{r}-\bm{r}_h\|_0$} &
		\multicolumn{2}{c|}{$\|\bm{u}-\bm{u}_h\|_0$} &
		
		\multicolumn{2}{c| }{$\|p-p_{h}\|_0$} &
		\multirow{2}{*}{DOF} \\
		\cline{3-8}
		
		&  &Error &Rate  &Error &Rate  &Error &Rate  \\
		\hline
		&2	&2.06E-01	&	    &3.50E-01	&	    &1.76E-01	&	    	&846\\
		&4	&1.40E-01	&0.55 	&2.75E-01	&0.35 	&1.13E-01	&0.63 		&5976\\
		$\bm{\mathcal{P}}^{{\color{black}\rm div}}_{\Gamma,k}\bm{g}_T$&8	&7.55E-02	&0.90 	&2.10E-01	&0.39 	&5.59E-02	&1.02 		&44640\\
		&16	&3.80E-02	&0.99 	&1.51E-01	&0.48 	&2.35E-02	&1.25 		&344448\\
		&32	&2.50E-02	&0.60 	&1.05E-01	&0.53 	&9.50E-03	&1.30 		&2704896\\
		
		\cline{1-9}
		&2	&2.62E-01	&	    &2.92E-01	&	    &1.62E-01	& 	    	&846\\
		&4	&2.33E-01	&0.17 	&2.43E-01	&0.27 	&1.10E-01	&0.56 		&5976\\
		$\bm{\Pi}^{\partial}_k\bm{g}_T$&8	&2.42E-01	&-0.06 	&1.92E-01	&0.34 	&5.43E-02	&1.02 		&44640\\
		&16	&3.02E-01	&-0.32 	&1.40E-01	&0.46 	&2.29E-02	&1.24 		&344448\\
		&32	&4.00E-01	&-0.40 	&9.78E-02	&0.52 	&9.30E-03	&1.30 		&2704896\\
		\Xhline{1pt}
	\end{tabular}
\end{table}

We conclude that all numerical experiments in this section verify that the underlying HDG scheme is appropriate for solving the Maxwell model problem \eqref{source} with both high and low regularities. 

\color{black}{
	\section{Conclusion}
	We have proposed and analyzed a new HDG method for solving the Maxwell operator. Theoretical results are obtained under a more general regularity requirement. In particular, a special treatment is applied to approximate data on the boundary for the low regularity case.
	
	Our future work includes the study of fast solvers for the proposed HDG method for the Maxwell operator.
	We will also consider the application of analysis procedure presented in this work to the numerical solution of incompressible magnetohydrodynamics \cite{Xu-submit}.

	%Our future work includes the study of fast solvers for the proposed HDG method in this paper. Moreover, we would like to develop an accurate, efficient and robust HDG method for solving time-harmonic Maxwell equations with high frequency. For the high frequency (large wavenumber) case, we will consider the derivation of the wavenumber-explicit regularity results. We also wish to establish the stability and error estimates for HDG method, where the coefficient constants are desired to be independent of the wavenumber.
} \color{black}

\bibliographystyle{siam}
\bibliography{mybib}{}

\end{document}